\documentclass[journal]{IEEEtran}
\usepackage{algorithmic}
\usepackage[ruled,vlined]{algorithm2e}
\usepackage{xcolor}
\usepackage{caption}
\usepackage[utf8]{inputenc}
\usepackage{amsmath}
\usepackage[colorinlistoftodos]{todonotes}
\usepackage{multirow}
\usepackage{graphicx}
\usepackage{amsmath}
\usepackage{breqn}
\usepackage{subcaption}
\usepackage{float}
\usepackage{enumerate}
\restylefloat{table}
\usepackage[utf8]{inputenc}
\usepackage{enumitem}
\usepackage{verbatim}
\usepackage[font=small,labelfont=bf]{caption}
\usepackage{tabularx}
\usepackage{makecell}
\usepackage{bbm}
\usepackage{cleveref}

\title{Public Transit for Special Events: \\
Ridership Prediction and Train Optimization}
\author{Tejas Santanam, Anthony Trasatti, Pascal Van Hentenryck, and Hanyu Zhang 
\\ Georgia Institute of Technology, Atlanta
}
\date{\today}

\begin{document}

 \maketitle

 \begin{abstract}
Many special events, including sport games and concerts, often cause surges in demand and congestion for transit systems. Therefore, it is important for transit providers to understand their impact on disruptions, delays, and fare revenues. This paper proposes a suite of
data-driven techniques that exploit Automated Fare Collection (AFC) data for evaluating, anticipating, and managing the performance of
transit systems during recurring congestion peaks due to special events. This includes an extensive analysis of ridership of the
two major stadiums in downtown Atlanta using rail data from the Metropolitan Atlanta Rapid Transit Authority (MARTA). The paper first
highlights the ridership predictability at the aggregate level for each station on both event and non-event days. It then presents an
unsupervised machine-learning model to cluster passengers and identify which train they are boarding. The model makes it possible to evaluate system performance in terms of fundamental metrics such as the passenger load per train and the wait times of riders. The paper also presents linear regression and random forest models for predicting ridership that are used in combination with historical throughput analysis to forecast demand. Finally, simulations are performed that showcase the potential improvements to wait times and demand matching by leveraging proposed techniques to optimize train frequencies based on forecasted demand. 
 \end{abstract}

 





\section{Introduction}
\label{sec:intro}


Special events, including sports games, concerts, and
festivals, are important for transit providers; they often lead to
fundamentally different ridership patterns and bring significant fare revenues. In addition, special events may be the introduction of certain riders to a transit system, and hence it is critical to ensure that the system is smooth and efficient, the waiting times are reasonable, and the vehicles are not too crowded, in order to attract additional recurring ridership. 

This paper originated as a study of special events for the Metropolitan Atlanta Rapid Transit Authority (MARTA), the transit system of the city of Atlanta in the state of Georgia. In particular, the study aims at addressing two main objectives of MARTA:

\begin{enumerate}
\item Is is possible to forecast special event rail ridership based on expected attendance and the type of event using historical Automated Fare Collection (AFC) data;
\item Can optimizing train frequencies based on forecasted ridership significantly improve passenger wait times and congestion following the events?
\end{enumerate}

Tackling these overall objectives requires addressing some sub-objectives to the following questions about special event ridership: 
\begin{enumerate}[noitemsep]
    \item How many passengers are using the rail to travel to and from special events?
    \item Which stations are most popular for event riders?
    \item What are the passenger loads of the trains? 
    \item How long do riders wait for trains after an event? 
    \item Can special event ridership be predicted? 
    \item Does the incidence of nearby smaller events have an impact on ridership to larger special events?
\end{enumerate}

To answer these questions, the paper presents a suite of data-driven
methods that leverage Automated Fare Collection (AFC). It
demonstrates that the data-driven tools answer the above questions
on sporting events located at the two downtown stadiums, the
Mercedes-Benz Stadium and the State Farm Arena which host the Atlanta
Hawks, Atlanta Falcons, and Atlanta United. From a high level perspective, this paper proposes the following methodology for MARTA, and possibly other transit agencies:
\begin{enumerate}
    \item Obtain an attendance prediction from stadium ticket sales and predict ridership from that attendance using supervised machine learning;
    \item Using historical trends, estimate the arrival times of passengers to the stations near the stadium;
    \item Based on the forecasted arrivals, optimize the train frequencies to minimize waiting times, maximize safety, and minimize costs.
\end{enumerate}

To develop this methodology, the paper makes the following technical and analytical contributions.

\begin{enumerate}[noitemsep]

\item The paper introduces station signatures to highlight the
  ridership consistency at the aggregate level for the baseline
  (non-event) days and post-games event peaks. These consistencies 
  provide a foundation to estimate special-event ridership accurately 
  and answer the above questions.

\item The paper combines station signatures with origin-destination
  (OD) pairs to highlight the strong correlation between the origins
  of special-event ridership and parking locations.

\item The paper uses unsupervised learning and simulation to estimate
  the train loads, departure times, and the associated waiting times of the riders. These estimates were validated with a simulation model for train boardings, inferring the train occupancy.

\item The paper uses supervised learning to predict ridership with
  high accuracy, giving planners the tools to optimize train schedules
  for future events. The paper describes the performance of three
  predictive models: a linear regression, a random forest, and a
  combination of linear regression and random forest, where the random
  forest predicts the residual errors of the linear regression.

\item The paper demonstrates that optimizing the train frequencies
  based on the forecasted demand may improve wait times and significantly reduce the number of riders left behind during post-game peaks while using a similar number of trains. 

\end{enumerate}


The rest of the paper is organized as follows. Section
\ref{sec:LitReview} reviews prior work on similar topics. Section
\ref{sec:data} presents the case study. Section
\ref{sec:RidershipAnalysis} presents the analysis of the baseline
ridership and the signatures for entries and exits at rail
stations on weekdays and weekends. It also analyzes the special
event ridership and estimates how many riders use the rail to travel to and from the events and where they come from.  Section
\ref{sec:train-scheduling-capacity-analysis} applies unsupervised
learning to cluster riders and predict which train they board, their
waiting times, as well as the departure times of the trains. Section
\ref{sec:ridership-prediction} presents simple supervised learning
models to predict the ridership for various types of recurring special
events that can be used with historical trends to forecast arrivals at nearby stations after the game. This section demonstrates how optimizing the train frequencies based on this forecasted demand may significantly improve the post-game congestion and wait times. This section also addresses the
question of how nearby smaller events affect the ridership to larger
events.

\section{Literature Review}
\label{sec:LitReview}
Automated Fare Collection (AFC) technologies have enabled more
sophisticated analysis of transit ridership \cite{barry2002origin, 7558221}. Various data sources have been used to study special-event
ridership including survey, AFC, and web data \cite{doi103141224604, camara2012internet}.
Rodrigues et al. (2017) use a Bayesian additive model to understand
and predict event riders arriving by public transit using Singapore
smart card data \cite{7765036}. That model creates a separate
prediction for each different event and baseline ridership using 5
months of data, grouping arrivals in 30 minute bins. Karnberger et
al. (2020) look at the Munich public transit system, also using some
AFC data. They look at weekly system averages and build a gradient
boosted random forest prediction system for ridership between linked
stations in parts of Munich \cite{KARNBERGER2020102549}. The type of
day (holiday, weekend, etc.) and the existence of a few types of
events are used as inputs to the model which focuses on link-level
riders.

Other papers looked at prediction of sporting event
attendees. Ni et al. (2017) establish a correlation between tweets
related to an event and event ridership flow. The paper builds a
linear regression model to predict ridership from the number of tweets for Mets games and US Open tennis matches \cite{7583675}. King (2017) predicts NBA game attendance using random forest models, but the ridership and travel modes are not considered \cite{King2017PredictingNB}.

Short-term prediction has also been the focus of other ridership
models, which contrasts to the longer time horizon considered in this
paper. Li et al. (2017) use a multiscale radial basis function network
to predict rail ridership at three large Beijing rail
stations. The model predicts riders at a station half an hour in the
future using a one-step-ahead model \cite{LI2017306}.

\section{The Case Study}
\label{sec:data}
\subsection{Map of MARTA Rail System}

\begin{figure}[t!]
\centering
\includegraphics[width=.75\linewidth]{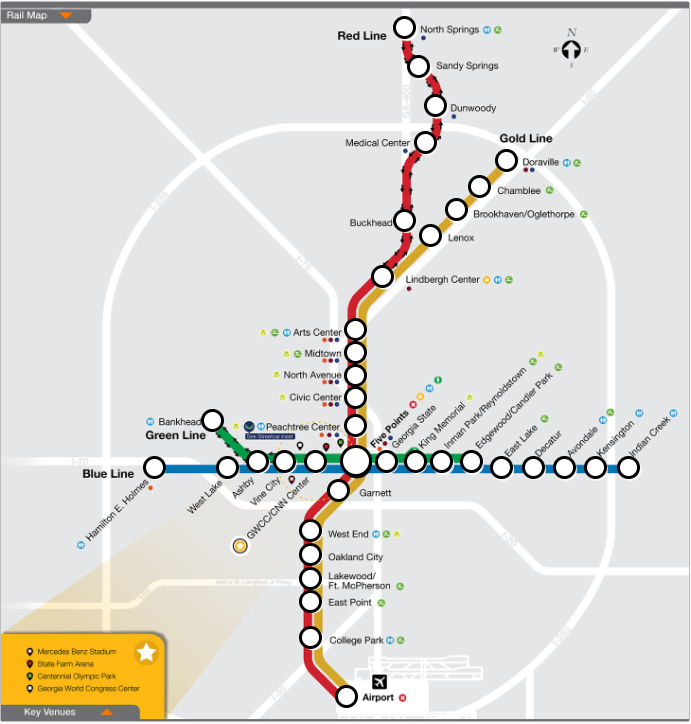}
\caption{The MARTA Rail Lines and Stations}
\label{fig:MARTAmap}
\end{figure}

Figure \ref{fig:MARTAmap} depicts the four MARTA rail lines \cite{MARTA-map}. The Red
and Gold lines run North-South and the Blue and Green lines run
East-West, with the two directions intersecting at Five Points. The
three event locations considered in this paper are
\begin{enumerate}[noitemsep]
    \item The Mercedes-Benz Stadium;
    \item The State Farm Arena;
    \item The Georgia World Congress Center.
\end{enumerate}

\noindent
These three key venues located in downtown Atlanta on the Blue and Green lines are highlighted in the left center of the map. The two closest stops are Dome/GWCC and Vine City. Users coming from the North or South can use the Red or Gold
line and transfer at Five Points to get on the Blue or Green line. These
locations are important for the subsequent analyses.

\subsection{Event Data}
\label{subsec:Eventdata}

The event data provided by MARTA is a list of many public
events in the Atlanta area during 2018 and 2019 containing
\begin{enumerate}[noitemsep]
    \item the type of Event;
    \item the location;
    \item the date;
    \item the time;
    \item the estimated attendance. 
\end{enumerate}
An example of event entries is presented in Table \ref{tab:Event Data Example}.

\begin{table}[!t]
\centering
\resizebox{\linewidth}{!}{%
\begin{tabular}{|c|c|c|c|c|}
\hline
Begin Date & Category & Event & Location & Attendance \\ \hline
\begin{tabular}[c]{@{}c@{}}01/15/2018 \\ 15:00:00\end{tabular} & Basketball - Hawks & Hawks v San Antonio Spurs & State Farm Arena & 15,000 \\ \hline
\begin{tabular}[c]{@{}c@{}}01/16/2018  \\  07:00:00\end{tabular} & Conference & Mary Kay Leadership Conference & \begin{tabular}[c]{@{}c@{}}Georgia World \\ Congress Center\end{tabular} & 7,000 \\ \hline
\begin{tabular}[c]{@{}c@{}}01/16/2018  \\  09:00:00\end{tabular} & AmericasMart & Intnl Gift \& Home Furnishings & AmericasMart & 72,000 \\ \hline
\end{tabular}%
}
\caption{Event Data Example}
\label{tab:Event Data Example}
\end{table}

The three most popular venues for large events in Atlanta are the
Georgia World Congress Center, the Mercedes Benz Stadium (MBS), and
the State Farm Arena, which were mentioned earlier as the focus of
this paper. In 2018 and 2019, there were 1330 special events with an
estimated attendance greater than 500 people in Atlanta. 706 of these
1330 special evemts were held in the Georgia World Congress Center,
the Mercedes Benz Stadium (MBS), or the State Farm Arena. The three
locations are geographically close to each other; moreover, the
closest two rail stations are the Dome/GWCC and Vine City stations.

This paper focuses on events with the largest impact and, in
particular, the 200 sporting events (Basketball, Soccer, and Football games).  All basketball games were held in the State Farm Arena, while
the soccer and football games were held in the MBS. Apart from the
sporting events, there were 74 conferences, 53 conventions, and 102
expos \& shows at the target locations. However, these events
generally have smaller attendances, and riders leave these events in
patterns that are fundamentally different from sporting events. In
many cases, people go in and out of the events throughout the day,
leading to a  ridership more dissipated over a larger time
horizon. Overall, these events have a smaller impact on the congestion of the rail system. These factors, as well as low sample sizes, make longer-term events like conventions and conferences poor subjects for analysis. The event data is summarized in Table
\ref{tab:dataOverviewNew}.

\begin{table}[!t]
\centering
\resizebox{\linewidth}{!}{%
\begin{tabular}{cccc}
\hline
\multicolumn{1}{|c}{Primary Event} & \# of Days &
\begin{tabular}[c]{@{}c@{}}Avg. Attendance\\ \end{tabular} & \multicolumn{1}{c|}{\begin{tabular}[c]{@{}c@{}}Avg. Post-Game\\ Ridership \end{tabular}} \\ \hline
Basketball -Hawks & 67 & 15,278 & 1,425 \\
Football Games & 28 & 69,477 & 10,846 \\
Soccer & 39 & 52,712 & 8,037 \\ \hline
\multicolumn{1}{|c}{Day Type} & \# of Days & \begin{tabular}[c]{@{}c@{}}Avg. Attendance\\ of Primary Event\end{tabular} & \multicolumn{1}{c|}{\begin{tabular}[c]{@{}c@{}}Avg. Post-Game\\ Ridership \end{tabular}} \\ \hline
Single Event & 78 & 38,062 & 5,487 \\
Two Events & 56 & 36,712 & 5,082
\end{tabular}%
}
\caption{The Event Data Overview.}
\label{tab:dataOverviewNew}
\end{table}

\subsection{Automated Fare Collection Data}

To enter or exit the MARTA rail system, customers are required to use
a ticket or a reloadable card (the ``Breeze Card'') at the gates of
individual stations. MARTA provided anonymized transaction-level data
showing tap-in and tap-out times and locations for the rail network
from 2016 to mid-2020. Example entries of Breeze Card data are shown
in Table \ref{tab:RailDataExample}.

\begin{table}[!t]
\centering
\resizebox{\linewidth}{!}{%
\begin{tabular}{|l|l|l|l|}
\hline
\textbf{UserId} & \textbf{TransactionDT} & \textbf{UseType} & \textbf{Station} \\ \hline
101	& 2018/3/20 8:01 &	Entry (Tag On) &	Doraville \\ \hline
101	& 2018/3/20 8:28 &	Exit (Tag Off) &	Lindbergh Center \\ \hline
101	& 2018/3/20 20:05 &	Entry (Tag On) &	Lindbergh Center \\ \hline
101	& 2018/3/20 20:21 &	Exit (Tag Off) &	Doraville \\ \hline
\end{tabular}%
}
\caption{An example of AFC rail data.}
\label{tab:RailDataExample}
\end{table}

\begin{table}[!t]
\centering
\resizebox{\linewidth}{!}{%
\begin{tabular}{|l|l|l|l|l|}
\hline
\textbf{UserId} & \textbf{EntryDT} & \textbf{EntryStation} & \textbf{ExitDT} & \textbf{ExitStation} \\ \hline
101	& 3/20/18 8:01 &	Doraville & 3/20/18 8:28 &	Lindbergh  \\ \hline
101	& 3/20/18 20:05 &	Lindbergh  & 3/20/18 20:21 &	Doraville \\ \hline
\end{tabular}%
}
\caption{An example of chained AFC Rail data.}
\label{tab:ChainedRailDataExample}
\end{table}

Trip chaining is performed to turn these individual transactions into Origin-Destination (OD) pairs. Tap-ins and tap-outs, when chained, are
sufficient to determine where a rider enters and exits the rail
network. For the most part, entries are matched to the following exit and create an OD pair. Table \ref{tab:ChainedRailDataExample} shows an example of the chained trips after they are processed. In some cases, the chained trip is an entry and exit at the same station back to back. For example, when an exit tap occurs, but there is a missing entry tap, then the system records a ``forced entry" transaction and an exit transaction at the exit station at the same time. Similarly, if someone tries to re-enter a station, but there is a missing exit transaction or an extended amount of period has elapsed (3-4 hours), the system will add a ``forced exit" transaction at the most recent tap-in location. Since the majority of transactions follow the expected pattern (pairs of distinct locations), the analysis focuses on these transactions.

\section{Ridership Analysis}
\label{sec:RidershipAnalysis}

This section first analyzes ridership on baseline (non-event) days, as
average day patterns help identify the effects of special events on
the system. The analysis considers baseline days and creates
individual ridership signatures for each station.  These signatures
are then used to calculate the ridership that can be attributed to
special events. The analysis also assesses the consistency of
special-event ridership.


\subsection{Station Signatures}

To create representative signatures for baseline days, the transaction
data is to partition on the basis of
\begin{itemize}[noitemsep]
\item weekday versus weekend;
\item entry versus exit
\item event day versus baseline day.
\end{itemize}

\noindent This partitioned data is used to create four baseline
signatures for each station: weekday entry, weekday exit, weekend
entry, and weekend exit. Since the rail is closed between 1:30am and
4:30am, the day is defined as the 24-hour period starting at 3AM and
finishing 3AM to capture riders returning after midnight. Each day is
further partitioned into 15 minute intervals and the number of riders
are counted for each interval for each of the four types of
transactions.




\begin{figure}[t]
\centering
\includegraphics[width=.7\linewidth]{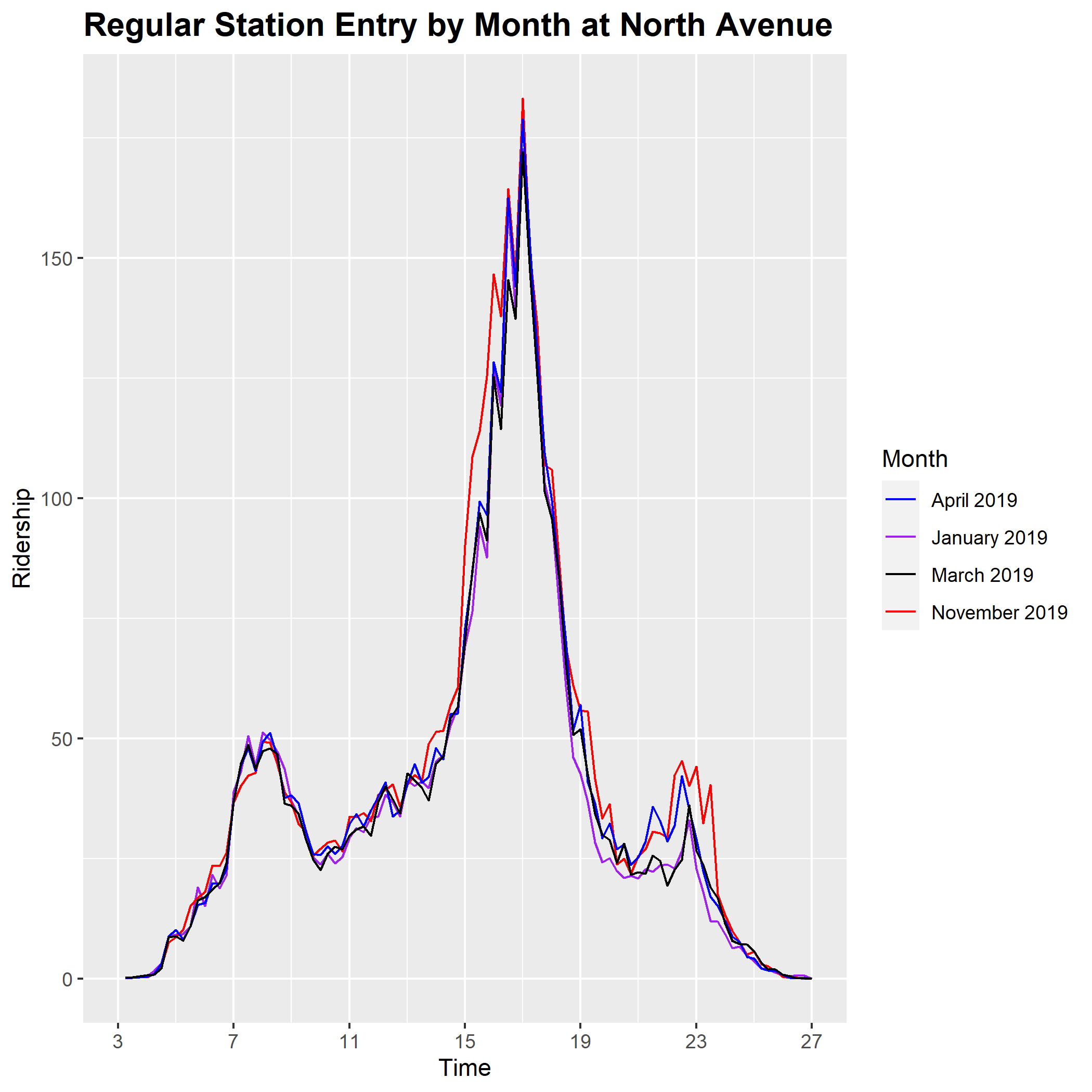}
\caption{Monthly Baseline Signatures for Entries at North Ave Rail Station}
\label{fig:monthconsistency}
\end{figure}

Figure \ref{fig:monthconsistency} shows how the signatures are very similar and consistent month-to-month. This kind of consistency is good for future modeling and planning.
In the appendix, Figures \ref{fig:BaselineStart}--\ref{fig:BaselineEnd} depict the four
baseline signatures for the North Avenue and Midtown stations. The shaded region represents the 10th-90th percentile range for each bin. The
signatures show regular commute spikes on weekdays for both entry and
exit signatures. They also show that riders follow similar patterns
throughout the year. Figures
\ref{fig:BaselineSignatureVineCityEntry}--\ref{fig:BaselineSignatureVineCityExit}
report the results for the Vine City station, which does not have many
regular commuters. There is still a weekday and weekend difference for
the ridership of Vine City station, but the overall magnitude of
ridership is low compared to the North Avenue and Midtown stations.
As will become clear, this low ridership changes drastically in
presence of an event as it is one of two closest stations to the
nearby Mercedes-Benz Stadium and State Farm Arena.

\subsection{Event Ridership Estimation} 

The last section highlights the consistency of the baseline ridership. Pereira et al. measures and detects hotspots by counting the number of riders greater than the median, where the curve exceeds the baseline's 90th percentile \cite{7021960}. This section uses this same technique to measure the special event ridership. 
This analysis focuses on riders with an
origin or destination at the Dome/GWCC and Vine City stations, since the majority of riders use these two stations to get to and from the events at the three venues of focus. 
This measurement is formalized as follows. Let

\begin{itemize}[noitemsep]
    \item $r_b(t)$ be the baseline ridership at time $t$;
    \item $r^+_b(t)$ be the 90th percentile of the baseline;
    \item $r_e(t)$ be the actual ridership during the event at time $t$;
    \item $t^s$ be the earliest time with $r_e(t) > r^+_b(t)$ to capture the start of the event ridership;
    \item $t^e$ be the latest time with $r_e(t) > r^+_b(t)$ to capture the end of the event ridership;
    \item $r_a(t)$ be the number of rail riders at time $t$ who attended the event;
\end{itemize}
The number of rail riders who attend the event at time $t$ is given by
\begin{equation}
r_a(t) = r_e(t)-r_b(t) \mathbbm{1}{(r_e(t)>r^+_b(t))}
\end{equation}
and the number of rail riders attending the  event is given by
\begin{equation}
T_a = \sum_{\forall t \in [t^s,t^e]}{r_a(t)}.
\end{equation}



\begin{figure}[!t]
 \centering
 \includegraphics[width=0.7\columnwidth]{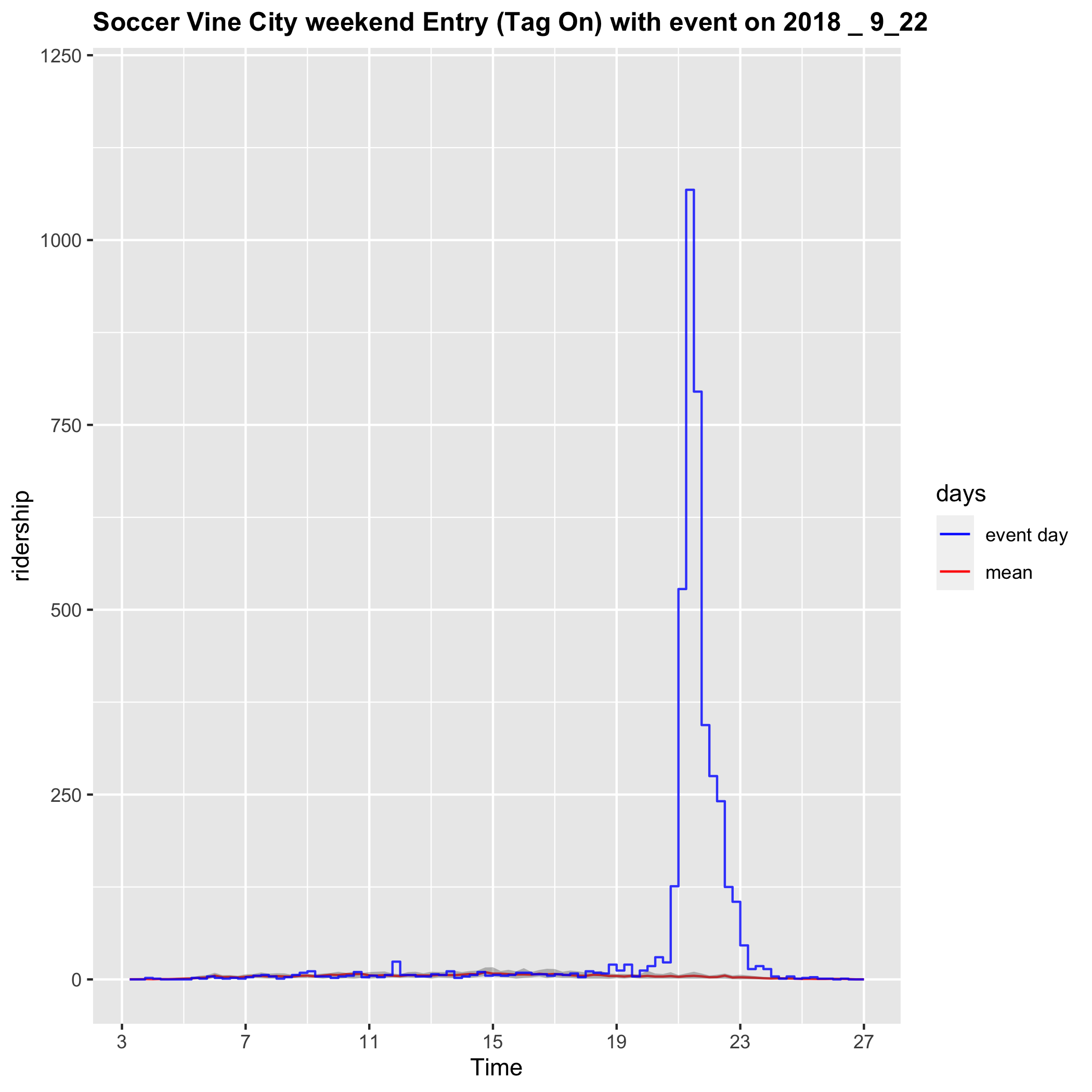}
 \caption{Vine City station's post-game ridership on September 22nd, 2018 versus its baseline weekend signature.}
 \label{fig:VineCity_PostGameRidershipSignature09_22_2018}
\end{figure}

\begin{figure}[!t]
 \centering
 \includegraphics[width=0.7\columnwidth]{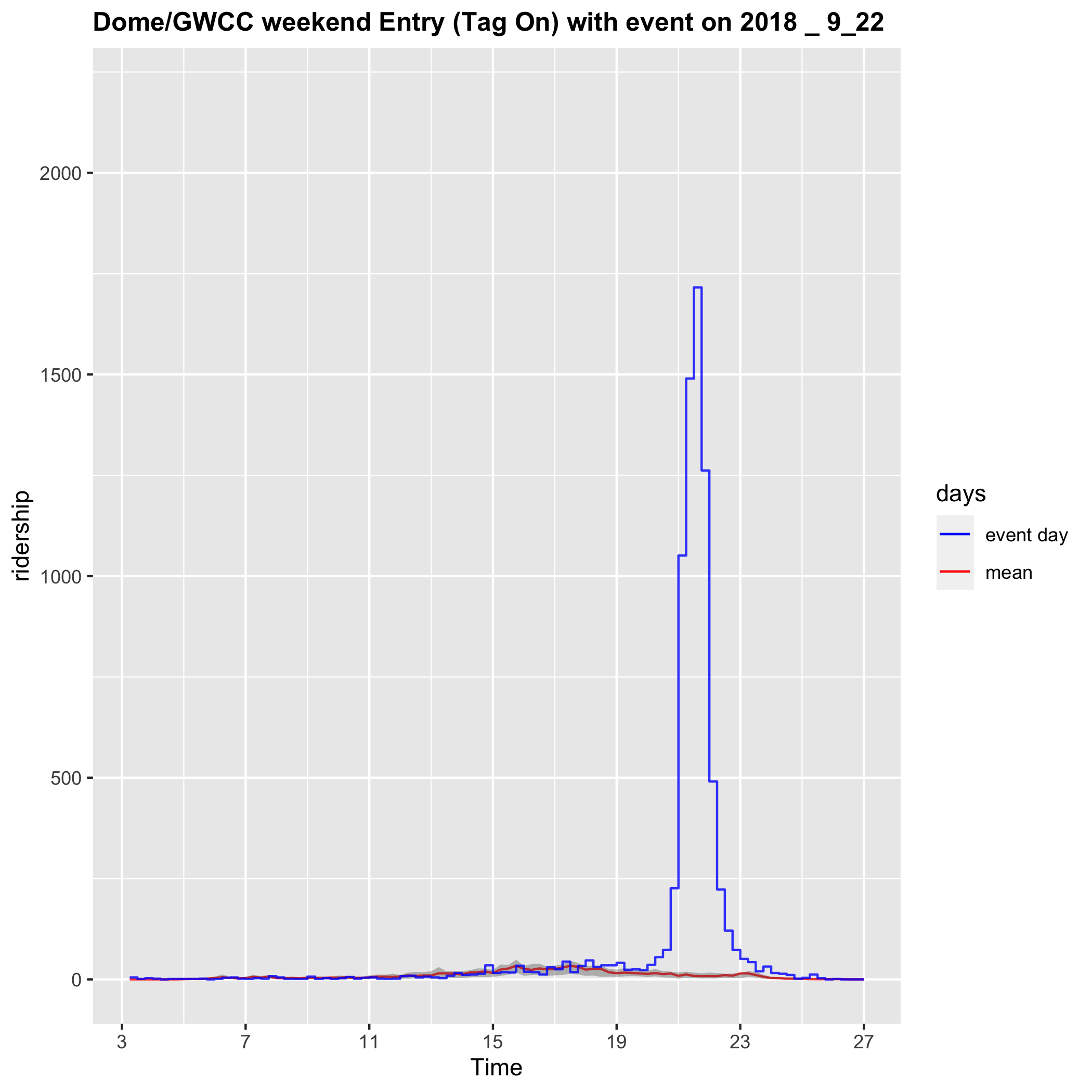}
\caption{Dome/GWCC station's post-game ridership on September 22nd, 2018 versus its baseline weekend signature.}
\label{fig:Dome_PostGameRidershipSignature09_22_2018}
\end{figure}



To illustrate these equations, consider the Atlanta United game on
September 22, 2018. Approximately 10,813 people entered
the Dome/GWCC and Vine City stations after the game, including
3,787 from Vine City and 7,026 from Dome/GWCC. The event had no
significant influence on the ridership on any other station. 
The two station signatures for September 22nd, 2018 
are presented in Figures
\ref{fig:VineCity_PostGameRidershipSignature09_22_2018} and
\ref{fig:Dome_PostGameRidershipSignature09_22_2018}.

\begin{figure}[!t]
\centering
\includegraphics[width=0.7\columnwidth]{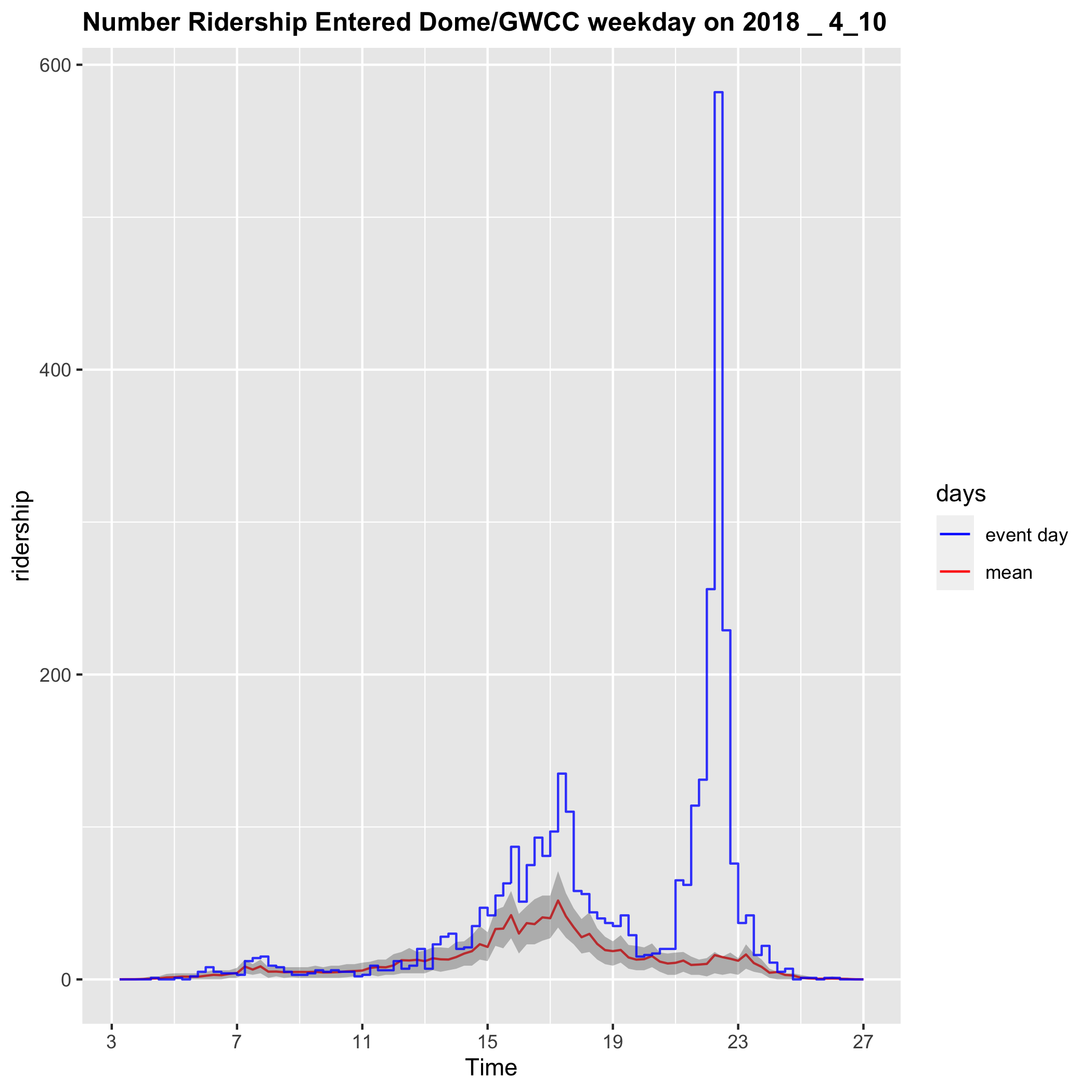}
\caption{Illustrating a Double Event Day.}
\label{Double Event Ex.}
\end{figure}

A double event day is shown in Figure \ref{Double Event Ex.}. 
The later spike corresponds to Hawks
v.s. Philadelphia 76ers game in State Farm Arena at 17:30. The earlier
spike corresponds to MODEX 2018, an expo in Georgia World Congress
Center, which starts at 10:00. It can be observed that the event departure pattern of
the riders for the non-sporting event depart is much more
spread-out than for sporting events.

The obtained values represent estimates for special-event ridership that are later used to build prediction models. It is assumed that the vast majority of these riders are indeed traffic due to the special events, as the ridership show large deviations above the normally low baselines at these stations.

\subsection{Event OD Patterns}
\label{sec:od-analysis}



This section focuses on the OD analysis of all Atlanta United games to understand which areas the special event ridership is coming from. The analysis can be used to understand the distributions for the origins (before the game) and the destinations (after the game) of these riders. Understanding where special-event riders come from can help transit agencies improve their offerings. The destination analysis suggests where riders might live, what forms of transportation they take, and what other factors contribute to the stations they use. Analysis of other types of games give similar results with changes mainly to the magnitude of station ridership.

\paragraph{Data} Without loss of generality, this section focuses on the destinations after Atlanta United games, since the later sections will focus on post-game service analysis and simulation. This post-game analysis focuses on riders who enter Dome/GWCC or Vine City stations 1 to 4 hours after the start of the game. 


\begin{figure}[!t]
    \centering
    \includegraphics[width=1\linewidth]{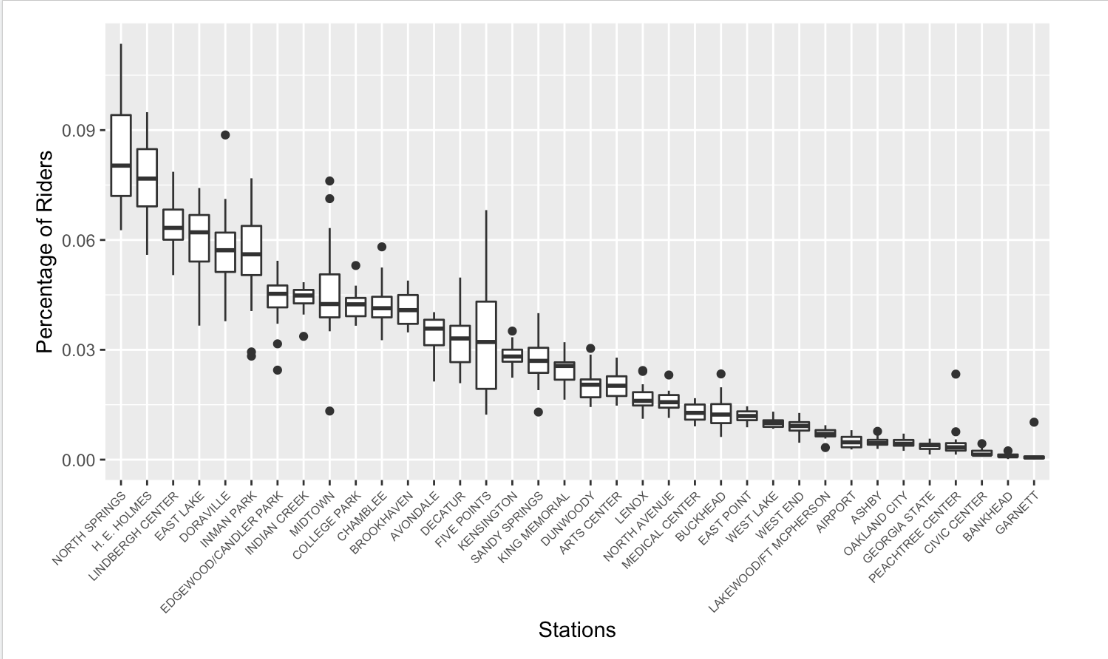}
    \caption{The Box Plot Showing the Median Percentages of the Post-Game Destinations}
    \label{BoxPlotpost-gamedistribution}
\end{figure} 

\begin{figure}[!t]
    \centering
    \includegraphics[width=.48\linewidth]{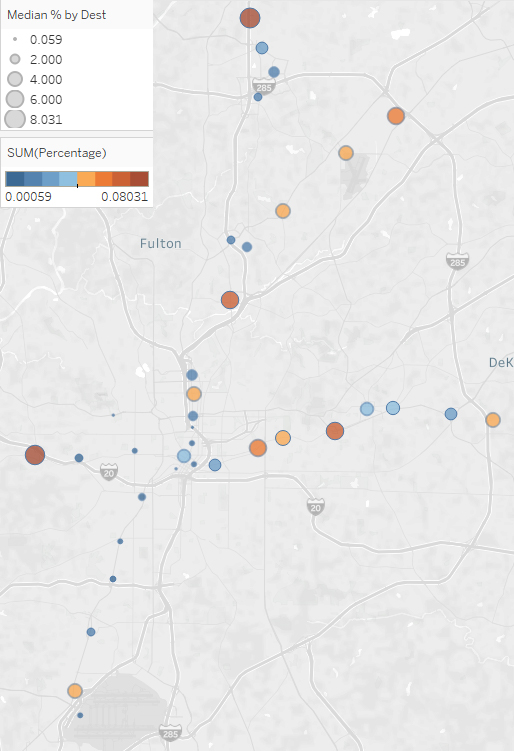}
    \caption{A heat map showing the median percentages for post-game destinations.}
    \label{Heatmappost-gamedistribution}
\end{figure}

Figure
\ref{BoxPlotpost-gamedistribution} shows the destinations of event
riders: the most frequent destinations are North Springs, H.E. Holmes,
and Lindbergh Center. For this analysis, the raw ridership counts for each game are normalized to obtain the percentages of riders alighting at the destination stations. For other events like football and basketball games, while the exact ranking of busiest stations may vary, the most used and least used stations are still roughly the same and handle similar percentages of riders. The stations with the largest baseline usage post-game also have the highest variances, largely a result of the exact attendances of games. The stations that get utilized the least have the lowest variances from not being used regardless of event attendance. Some of the variation could also stem from changes in baseline rider behavior if they are aware the game is going on. The Five Points station has a large variation for its magnitude, but that is a result of Five Points being a busy station that is the only transfer point on the entire MARTA network and could be influenced by happenings all over the system. These median percentage for each station are also shown on a map in Figure \ref{Heatmappost-gamedistribution} to give a geographical representation of this data.

\begin{figure}[!t]
    \centering
    \includegraphics[width=.48\linewidth]{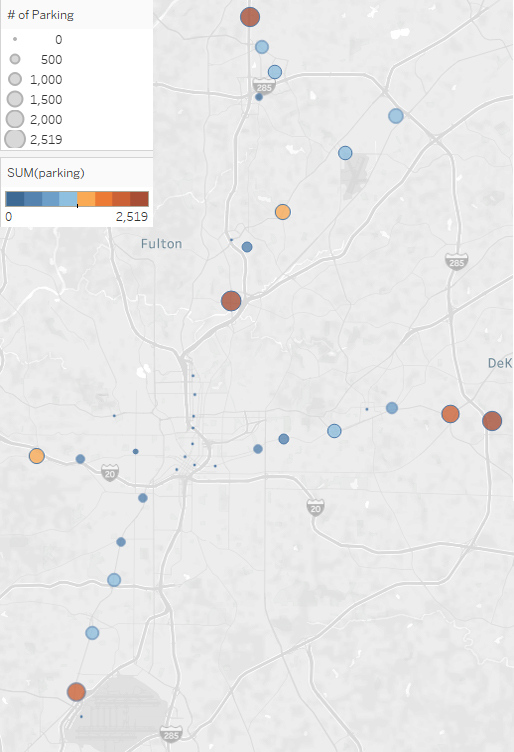}
    \caption{A heat map showing the number of parking spots per station.}
    \label{heatmapofparkingspots}
\end{figure}

Figure \ref{heatmapofparkingspots} shows a heat map
of number of parking spots available at each station. The heat map suggests that two major factors explain why riders use a particular station for events:
\begin{enumerate}[noitemsep]
    \item the proximity of the event location to the East/West line;
    \item the parking space availability.
\end{enumerate}
\noindent
The stations with a high number of riders are on the East/West line or
have ample parking or both. The four stations most used by riders
(i.e., North Springs station, H.E Holmes station, Lindbergh Center
Station, and East Lake Station) are all in the top six for parking
spots. A larger number of riders also use the East/West line although
the North/South line is closer to larger population centers.



\section{Entry-Exit AFC Analysis}
\label{sec:train-scheduling-capacity-analysis}

This section analyzes the performance of the train system after an
event. It shows how to estimate the train occupancy during the peak
post-event period when passengers are more likely to pack each
car. It also uses simulation to validate the estimated train capacity.

\subsection{Case Study Data}
\label{sec:Case Study Data}

The analysis in this section focuses on the Atlanta United game on September 22nd,
2018 for concreteness. Due to the nature of the rail system, the
primary focus is on riders using the rail to travel in the east
direction after leaving the stadium and entering either the Dome/GWCC
or the Vine City stations. Some passengers may proceed to switch to
another rail line, but the analysis focuses on the subset of
passengers who solely use the West to East tracks (either the Blue or
Green line).

\subsection{Train-Level Clustering}
\label{sec:TrainLevelClustering}

The train-level clustering analysis identifies the train schedules and
which riders were on the same train based on the AFC data. It consists
of three steps:
\begin{enumerate}

\item {\em Departure Time Inference:} The rider exit times 
  are used to obtain their departure times from the event stations;

\item {\em Rider Clustering:} Event riders are grouped in trains based
  on their inferred departure times;
  
\item {\em Schedule Estimation:} The train schedule after the game is
  estimated from the clustering results. The train departure time at a
  station is approximated as the latest arrival time of the riders on
  that train at that station.
\end{enumerate}

\subsubsection{Time Adjustment for Departure Time Inference}
\label{sec:Time Adjustment}

To determine the departure times of riders at the event stations,
their station exit times are shifted backwards, using the train travel
times. Given the data in Table \ref{tab:ChainedRailDataExample}, the EntryStation is denoted as $o$ and the ExitStation is denoted as station d. Let $e$ denote the event station, where station e is either station o or the trains are traveling from station o to station d through station e. A rider's travel time in the transit system can be decomposed by the following equation, where $\text{train Travel Time}_{od}$ is obtained from MARTA train schedule. 

\begin{equation*}
\begin{aligned}
    \text{ExitDT}_{d}= \text{EntryDT}_{o}+\text{Wait Time}_{o}\\
    +\text{Train Travel Time}_{oe}\\
    +\text{Train Travel Time}_{ed}
\end{aligned}
\end{equation*}

A departure time from the event station $e$ and a destination $d$ is defined as 
\begin{equation}
\begin{aligned}
        \text{Departure Time}_{e}=\text{ExitDT}_{d}-\text{Travel Time}_{ed}
\end{aligned}
\label{eq: departure time}
\end{equation}

As a event can have effect on multiple neighboring stations, to make the start time comparable. The arrival time at the event station $e$ from an origin station $e$ is defined as
\begin{equation}
    \begin{aligned}
        \text{Arrival Time}_{e}=\text{EntryDT}_{o}
        +\text{Travel Time}_{oe}
    \end{aligned}
    \label{eq: arrival time}
\end{equation}

\begin{figure}[!t]
    \centering
    \includegraphics[width=1\linewidth]{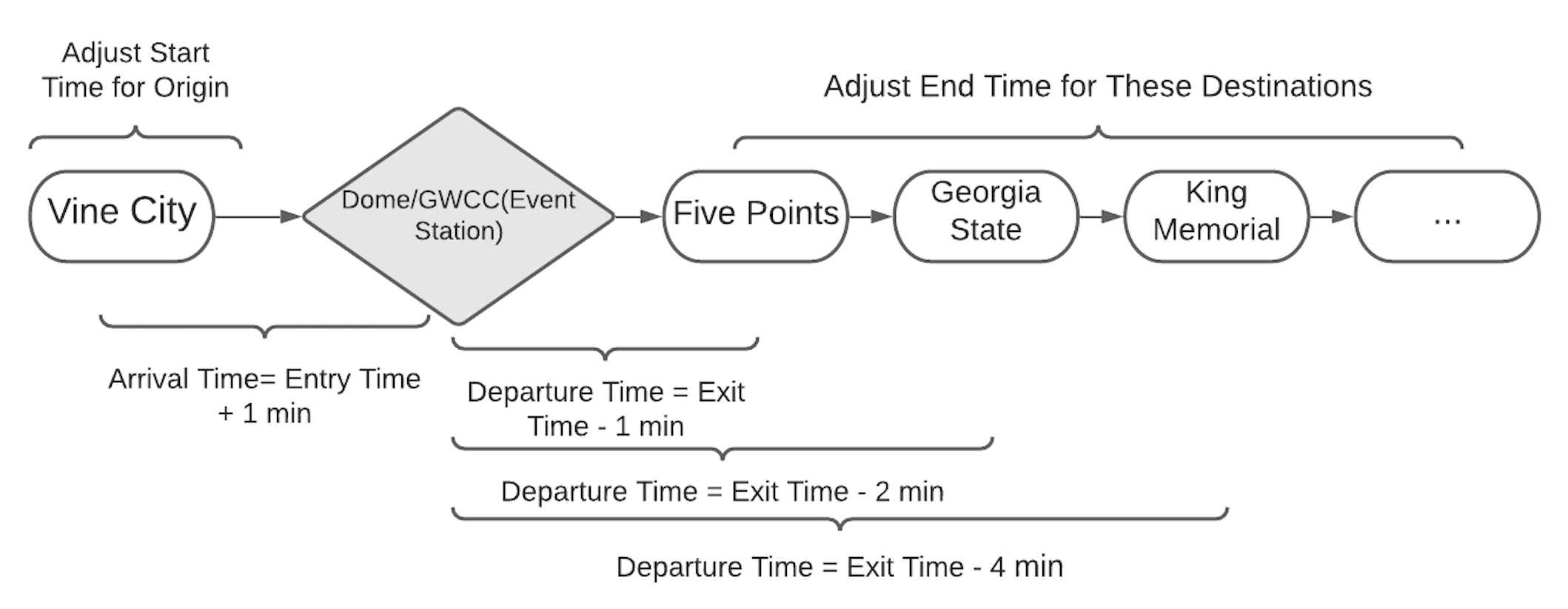}
    \caption{Departure Time Inference}
    \label{fig:Dep_time_inf}
\end{figure}

\begin{figure}[!t]
\centering
\includegraphics[width=1\linewidth]{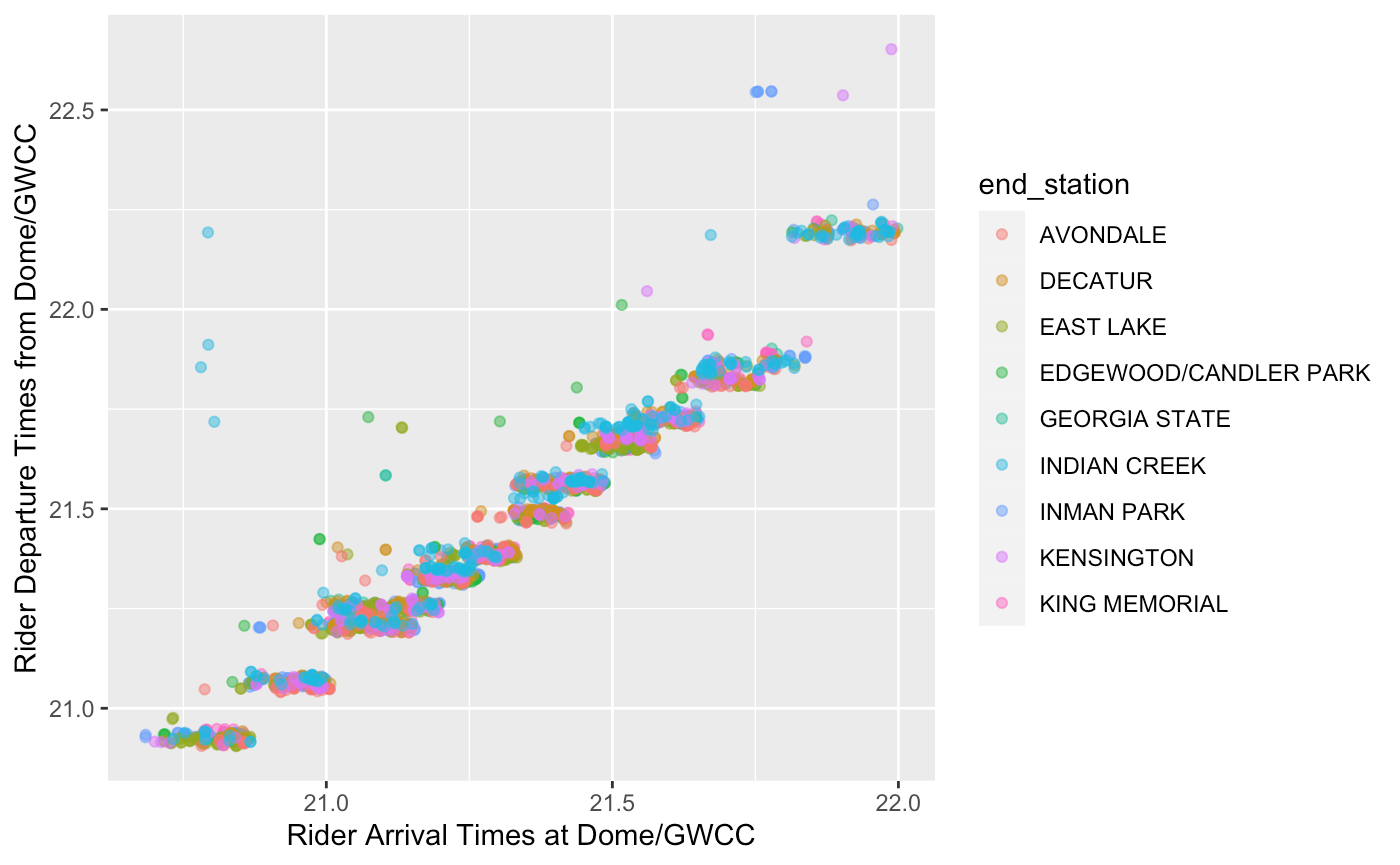}
\caption{Arrival and Departure Times of Riders at the Dome/GWCC Station.}
\label{fig:after-time-adjustment}
\end{figure}

Figure \ref{fig:Dep_time_inf} is an illustration of the departure time inference in this paper. The event station $e$ is Dome/GWCC station, and the riders boarding from Vine City station are also considered since the event also has an significant effect on the ridership Vine City. The riders' departure time and arrival time at the event station are adjusted according to equation \ref{eq: departure time} and \ref{eq: arrival time}. Figure \ref{fig:after-time-adjustment} shows departure time inference results with Dome/GWCC
as the event station. In the figure, the colors represent a different
alighting station for the riders. Observe the horizon clusters that
represent sets of riders boarding the same train.


\subsubsection{Rider Clustering}
\label{sec:Clustering}

Once the arrival and departure times are available, it is possible to
apply an unsupervised learning model to cluster riders in trains.
Algorithm HDBSCAN \cite{mcinnes2017hdbscan, malzer2019hybrid} was selected for this task,
because of its ability to obtain clusters of different densities
according to the mutual distances between the data points. This
section reports the results of this clustering for riders entering the
Vine City and Dome/GWCC stations after the game and alighting at
Edgewood/Candler Park, East Lake, Decatur, Avondale, and
Kensington. The most crowded period is between 20:40:00 and 22:00:00
and is the focus of this section.

\begin{figure}[!t]
\centering
\includegraphics[width=1\linewidth]{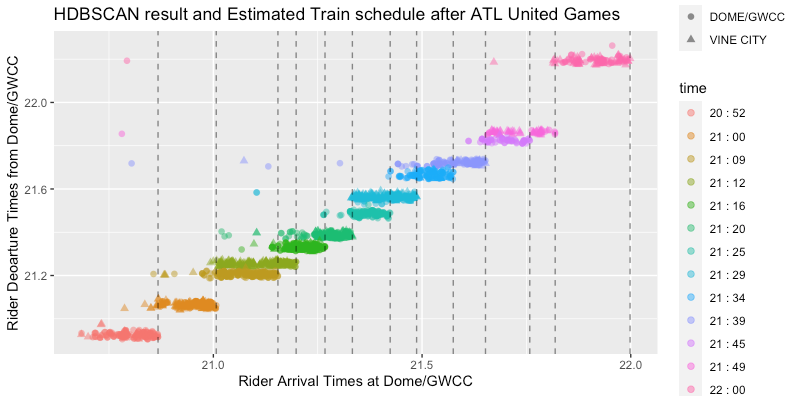}
\caption{Train Clusters for Riders Entering the Dome/GWCC and Vine City stations after the Atlanta United Game on September 22nd, 2018.}
\label{fig:20180922-clustering-result}
\end{figure}

To cluster the selected 2,392 riders in selected time interval,
HDBSCAN was run with its parameter {\tt MinPts} for minimum cluster
size set to $50$. An initial run detected 12 clusters and the 21:45
and 21:49 trains were not separated because the 21:45 train was
delayed at some stations and hence the departure time inference
resulted in some scattered data points. The first HDBSCAN run
identified the noisy data points and a second HDBSCAN run was applied
to produce the 13 clusters that correspond to the 13 trains that left
Dome/GWCC after the game.  The cluster results and the estimated train
arrival times at Dome/GWCC station are presented on Figure
\ref{fig:20180922-clustering-result}.  The estimated train departure
times at the Dome/GWCC station are plotted with dash lines which
represent the latest rider boarding time for each train.

The clustering algorithm assigns each rider to the corresponding train
they boarded. These passengers enter MARTA from two stations, Vine
City and Dome/GWCC. Because the trains are heading to Dome/GWCC from
Vine City, Vine City riders have a priority to board the train.
During peak times after the game, some riders cannot board the first
``available'' train. Moreover, some riders entering Dome/GWCC have to
wait up to three trains to board. Hence an important metric to
evaluate MARTA's performance is the percentage of riders left-behind
by each train, which is referred to as \textit{proportion-left-behind}
in this paper. Note that 3 of the 13 trains did not stop at the Vine
City station while all 13 trains did stop at Dome/GWCC station. Many
of the trains are mostly filled with Vine City passengers, so this
decision likely helps improve the wait times for riders using the
Dome/GWCC station.

\subsection{Train Capacity Utilization}
\label{sec:train-capacity-estimation}

To estimate how the trains are used, this section presents a simple
simulation model that estimates the percentage of riders left behind
at a station based on the arrival time of passengers and the train
capacity. By comparing the percentage of riders left behind computed
by the simulation model and the clustering, it becomes possible to
estimate how the trains are utilized. More precisely, the goal is to
find a train capacity that minimizes the distance between the outputs
of the simulation and clustering models, where the distance is
measured by mean absolute error loss.

Algorithm \ref{alg:left-behind-calculation} presents the simulation procedure:
It takes the train schedule estimated from the clustering model, the
arrival time of passengers, and a train capacity, and outputs an
estimation of the percentage of riders left behind by each train. For
instance, for the Atlanta United game, the rider arrival times are
broken into 13 intervals: riders arriving between 20:52 and 21:00 are
part of the demand for train at 21:00 and the simulation reports how
the percentage of riders who cannot board the next train.  The
algorithm uses the following notations:
\begin{itemize}[noitemsep]
    \item $I$ denote the number of trains;
    \item $i$ be the index of the trains;
    \item $S$ denote the set of all stations;
    \item $s$ be the index of the stations;
    \item $T_i$ denotes the departure time of train $i$;
    \item $C_i$ is the capacity of train $i$;
    \item $C_i^s$ is the remaining capacity of train $i$ at station $s$;
    \item $d_i^s$ is the number of passengers arriving at station $s$ between $T_{i-1}$ and $T_i$;
    \item $r_i^s$ is the number of riders who want to board train $i$;
    \item $l_i^s$ is the number of riders left behind by train $i$ at station $s$.
\end{itemize}

\begin{algorithm}[!t]
\SetAlgoLined
 $\text{left}_0^s$=0\;
  $C_i^0=C_i$\;
 \For{i \text{in} 1:I}{
 $r_i^s$ = $l_{i-1}^s$+ $d_i^s$ \;
 $l_{i}^s$= $\max(r_i-C_i^s,0)$\;
 $C_i^{s+1}=\max(0,C_i^s-r_i^s)$\;
 $\text{proportion-left-behind}_{i,s}=\frac{l_i^s}{r_i^s}$ \;
 }
\caption{Simulation}
\label{alg:left-behind-calculation}
\end{algorithm}

\begin{figure}[!t]
\centering
\includegraphics[width=1\linewidth]{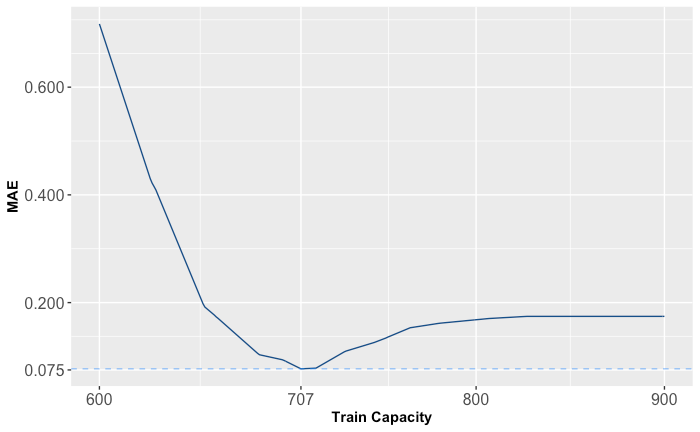}
\caption{Minimization of MAE Loss Function When Searching for the Maximum Capacity.}
\label{loss function result}
\end{figure}

\begin{table}[!t]
\centering
\resizebox{\linewidth}{!}{%
\begin{tabular}{|c|c|c|}
\hline
Train Time & \begin{tabular}[c]{@{}c@{}c@{}} Percentage left behind \\ estimated at Dome \\ by the clustering \end{tabular} & \begin{tabular}[c]{@{}c@{}} percentage left behind \\ at Dome in the simulation \\ 
(Maximum Capacity 707)\end{tabular} \\ 
\hline
20:52 & 0.04 & 0.00 \\ \hline
21:00 & 0.11 & 0.18 \\ \hline
21:09 & 0.08 & 0.19 \\ \hline
21:12 & 0.60 & 0.70 \\ \hline
21:16 & 0.29 & 0.29 \\ \hline
21:20 & 0.16 & 0.00 \\ \hline
21:25 & 0.09 & 0.00 \\ \hline
21:29 & 0.37 & 0.31 \\ \hline
21:34 & 0.35 & 0.17 \\ \hline
21:39 & 0.12 & 0.00 \\ \hline
21:45 & 0.07 & 0.00 \\ \hline
21:49 & 0.00 & 0.00 \\ \hline
22:00 & 0.00 & 0.00 \\ \hline
\end{tabular}%
}
\caption{Simulation Results for the Percentage of Riders Left Behind by each Train.}
\label{tab:simulation result}
\end{table}

Figure \ref{loss function result} depicts how the loss function
evolves for different train capacities, suggesting a train capacity
of 707. Table \ref{tab:simulation result} presents the results and
reports the percentages left behind by both the clustering and
simulation models for a train capacity of 707. The ``real''
percentages (from clustering) are larger than the simulated
percentages when the numbers are small. This is due to the fact that
the simulation assumes a perfect first-come-first-served
rule. However, this is not the case in real situation.
 
The maximum capacity of 707 is a lower bound estimation because riders
already on the trains (approximately a total of 35 people for the
whole time period) are not counted here. These trains are 6-car trains
post-game which have a recommended maximum capacity of 576
people. From this analysis, however, one can see that the maximum
capacity is often exceeded post-game: people often cram together in
very close quarters as a result. It is also likely that this
over-capacity situation leads to an increased risk of accident,
injury, or illness. However, the analysis simply confirms the
anecdotal evidence that people have a tendency to ``pack it in" after
sporting events. Note also that, under the assumption that people left behind end up
boarding the next train before the new arrivals, riders wait a maximum
of two trains, which corresponds to the case in Figure
\ref{fig:20180922-clustering-result}.




\section{Predictive Analytics}
\label{sec:ridership-prediction}

The end of a special event can lead to a large surge of passengers flooding to public transit. This can lead to crowding in the stations and extended wait times. The train operators adjust the normal schedule by pulling reserve trains ``out-of-pocket" in order to match the increased demand. Train operators have an important job of being aware of the state of the game, so that they can adjust the actual schedule accordingly. The goal of this section is to see if the proposed predictive analytics may be able to help train operators be better prepared for these recurring post-game demand spikes.

To optimize the train schedules of future events, it remains to demonstrate that the ridership can be predicted with high accuracy. This section shows that this is indeed possible. It first shows the consistency of station arrivals after a game. It then shows how to predict total ridership with high accuracy. Finally, these two results are leveraged as a demand forecast. The proposed train schedules designed using the forecasted demand are compared through simulation to the actual train schedules against the post-game arrival data. The proposed train schedules reduce the number of riders left behind and waiting times of the riders during post-game peaks.

\subsection{Throughput Consistency}

This section highlights the post-game throughput consistency for special events. Due to their larger sizes, events at Mercedes-Benz where the upper-deck seating was open are the focus of this section. In games with an open upper deck, an additional 30,000 seats are available for purchase in the upper deck of Mercedes-Benz Stadium. Due to the additional attendance, these events have a much larger impact on the MARTA rail system, especially compared to Atlanta Hawks games where the average attendance is only ~15,000 people. In this section, the post-game throughput is analyzed from 40 minutes before the end time to 80 minutes after the end time.

  In most games, it is assumed that the end time is the average game length after the scheduled start time: 1 hour \& 50 minutes for soccer and 3 hours \& 10 minutes for football. However, a few of the end times were adjusted in this analysis because it was believed that the end time might have been delayed due to injuries, delayed starts, or overtime. The delay in the actual end time of the game compared to the end time calculated using the average game length is referred to as the \textit{offset}. For each game, the \textit{offset} is estimated by comparing the throughput curves in cases where there was a clear delay to the peak of the throughput. When there are delays to the game, such as overtime for a Falcons game, the train operator waits to make the necessary adjustments to the actual schedule.

\begin{figure}[!t]
        \centering
        \includegraphics[width=0.9\linewidth]{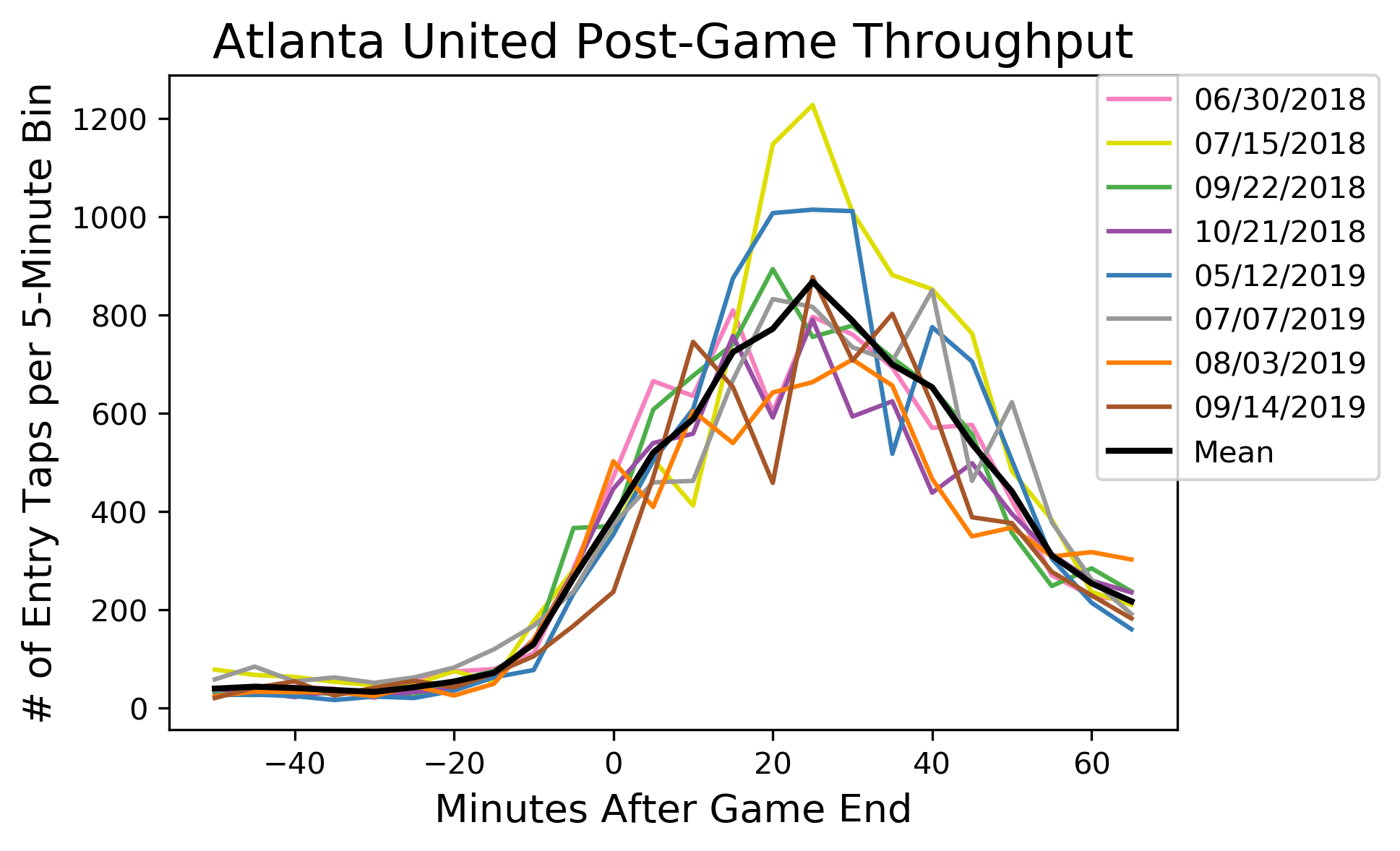}
        \caption{Post-game rider throughput at Dome/GWCC after Atlanta United games with the upper-deck seating opened.}
        \label{fig:atl-united-throughput}
\end{figure}

\begin{figure}[!t]
    \centering
    \includegraphics[width=0.9\linewidth]{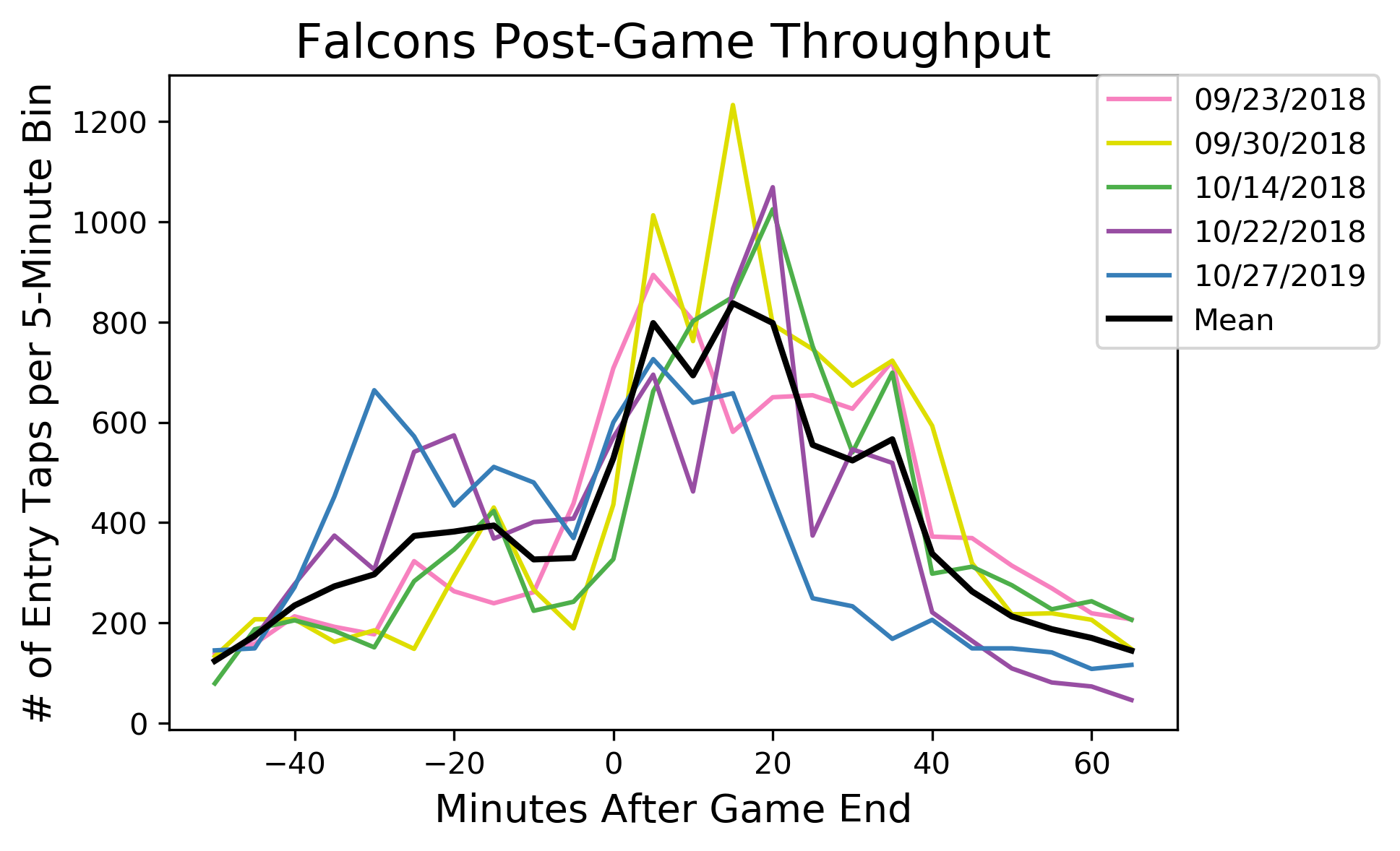}
    \caption{Post-game rider throughput at Dome/GWCC after Falcons games.}
    \label{fig:falcons-throughput}
\end{figure}

  Figure \ref{fig:atl-united-throughput} shows the entries to Dome/GWCC and Vine City are grouped into 5 minute bins and plotted for analysis. Note that three games had offset adjustments as stated later in Table \ref{tab:AtlUnitedSimulationTable}. This highlights the consistency of arrivals to the rail stations Dome/GWCC and Vince City rail stations after Atlanta United games with an open upper deck. The highest number of riders arriving in any bin is almost 1,200, which can be served with less than two trains assuming a train capacity of 707 as estimated in Section \ref{sec:train-capacity-estimation}. Note also that only ~8\% of people take a train going west, while the rest of the riders wait for trains going east. Figure \ref{fig:falcons-throughput} show a similar, yet distinct consistency for the football games. In some of the Falcons' games, there are some peaks that could align with the end of the third quarter (ex. 10/27/2019 where the falcons were down 24-0 at half time). Falcons games are typically more than an hour longer than Atlanta United games in length, so this could also explain why more fans leave early for Falcons games than for Atlanta United games.

The mean curve gives an relatively accurate estimation of the arrival patterns at the station. The mean curve is converted into percentages by dividing by the average ridership and used later in combination with the predicted ridership to obtain a predicted throughput curve. Given a ridership prediction, it is possible to estimate the arrival distribution over time at the station quite accurately.

\subsection{Ridership Prediction Model}
\label{sec:ridership-prediction-sub}

It remains to show that it is possible to predict the event ridership from the event attendance. Note that, for future events, a prediction of event attendance can be used in place of actual attendance. For example, NBA game attendance have been predicted using random forest models with a 6\% MAPE using team/opponent statistics, stadium capacity, local average income, team popularity and other factors in \cite{King2017PredictingNB}. The predictive model receives as inputs the attributes listed in Table \ref{tab:attributes}
and outputs $T_a$, i.e., the total number of riders for the event at
the considered station. A mean arrival curve (such as the
one shown in Figure \ref{fig:atl-united-throughput}) can then be used to
obtain an estimate of $r_a(t)$, i.e., the number of riders arriving at
time $t$. In the table, Event 1 is the event whose ridership must be predicted and Event 2 is
another event on the same day. For days with a single event, the
attributes of Event 2 are set to null. 
 
\begin{table}[!t]
\centering
\resizebox{\linewidth}{!}{%
\begin{tabular}{|c|c|p{3in}|}
\hline
 & \textbf{Attributes} & \textbf{Type} \\ \hline
\multirow{3}{*}{Event 1} & Category & Factor (Soccer/Football Game/Basketball) \\ \cline{2-3} 
 & Location & Factor (State Farm Arena/MBS/MBS\_Upper deck Open) \\ \cline{2-3} 
 & Attendance & Numeric, Attendance for Event 1 \\ \hline
 & wpdiff & Home Team win percentage minus Away Team win percentage \\ \hline
 & regularized margin & Margin of victory (loss) divided by standard deviation margin of victory for that league \\ \hline
\multirow{3}{*}{Event 2} & Category 2 & Factor (15 Categories) \\ \cline{2-3} 
 & Location 2 & Factor (GWCC/MBS/SFA/No  Location) \\ \cline{2-3} 
 & Attendance 2 & Numeric Attendance for Event 2) \\ \hline
\multirow{4}{*}{} & time\_difference & time distance (in minutes) of the two events, 0 if there is no second event \\ \cline{2-3} 
 & two\_event & Binary (True, if there is a second event) \\ \cline{2-3} 
 & week & Binary (True, if the day is weekend) \\ \cline{2-3} 
 & month & Factor (Month of the event) \\ \hline
\end{tabular}
}
\caption{The Input Attributes for the Predictive Models.}
\label{tab:attributes}
\end{table}


    
    
The results presented in this section focus on sporting events near
Dome/GWCC station and Vine City station, i.e., Atlanta Hawks games,
Atlanta Falcons games, and Atlanta United games, as there are enough
data points for those events to build strong models. These events are
also among those that had the largest impact on ridership. Post-game ridership is estimated using equation
$(1)$. The training data consists of 134 event days ranging from January 2018 to December 2019.

\begin{figure}[!t]
    \centering
    \includegraphics[width=\linewidth]{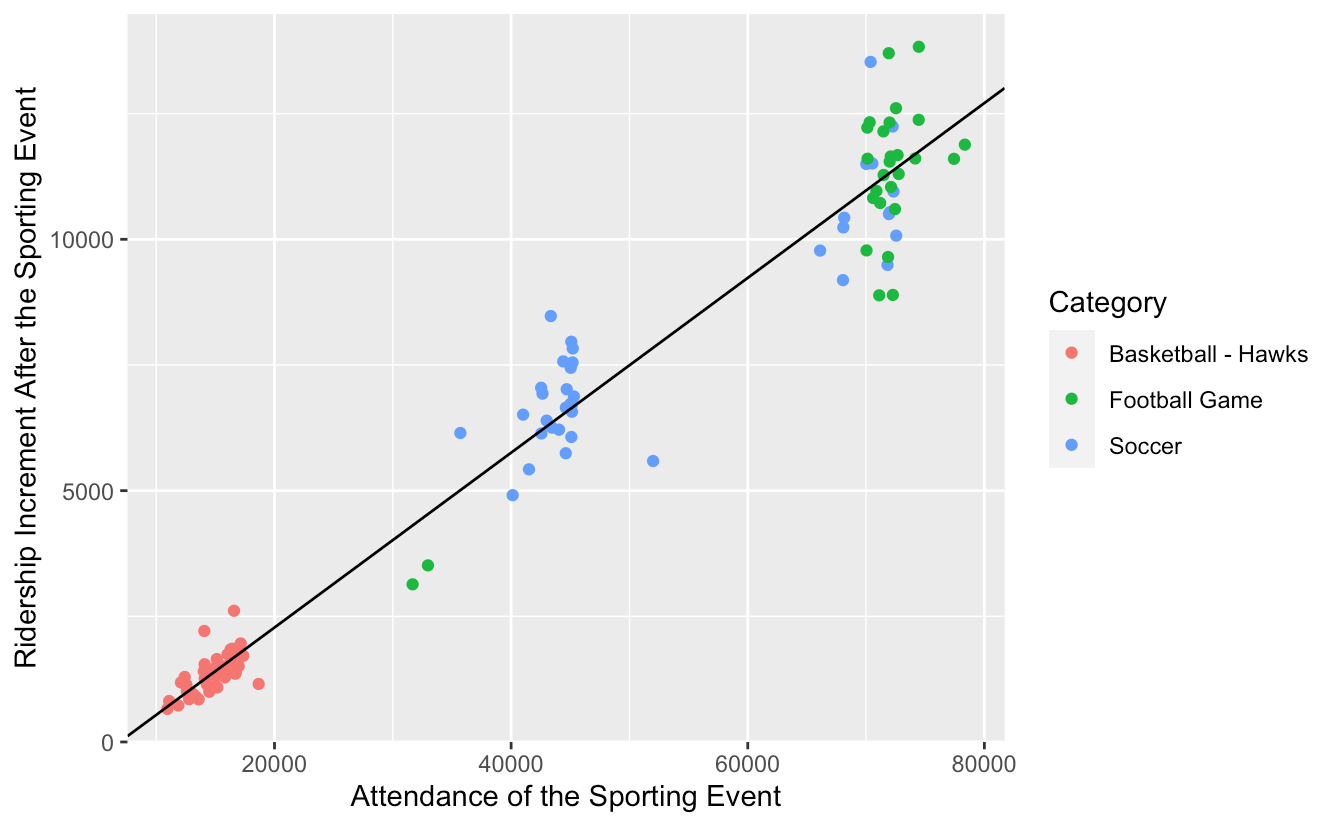}
    \caption{Linear trend between attendance and ridership increment after the main event}
    \label{fig:ML model linear trend}
\end{figure}

The first predictive model is a simple Linear Regression (LR)
\begin{equation}
\text{Ridership}=\beta_0+\beta_1\times \text{Attendance}
\end{equation}
that captures the strong linear dependence between strong linear
between the event ridership and the event attendance, which is
highlighted in Figure \ref{fig:ML model linear trend}.  The second
predictive model is a Random Forest (RF) that uses all the attributes
in Table \ref{tab:attributes}. The third model (LR+RF) is a
combination of the two: its goal is to fit a model
\begin{equation}
\text{Ridership}=-1201+0.1739\times \text{Attendance} + \epsilon
\end{equation}
where the residual $\epsilon$ is predicted by a random forest:
\begin{equation}
    \epsilon = \frac{1}{B}\sum_{b=1}^B T_b({\bf x})
\end{equation}
where $B$ is the number of decision trees, $T_b$ is the $b^{th}$
decision tree, and ${\bf x}$ is the input vector. $B$ is an
hyper-parameter obtained by fitting the model over different values
and selecting the one minimizing the RMSE. The LR+RF model recognizes
that the error term of the linear model depends on other factors,
e.g., whether there is a second event at that day or the win
percentage of the home team.

Table \ref{tab:result1} presents the results obtained using a {\em
  leave-one-out cross validation} because of limited sample size. RF
uses 1,500 trees and LR+RF uses $B=800$ trees. The proposed LR+RF
outperforms the other models among all metrics, i.e., MAE, MAPE, and
RMSE.

\begin{table}[!t]
\centering
\resizebox{0.8\linewidth}{!}{%
\begin{tabular}{|c|c|c|c|}
\hline
    & \textbf{MAE} &  \textbf{MAPE} & \textbf{RMSE} \\ \hline
    \textbf{LR} 
        & 509 
        & 0.1169    
        & 509.26    
    \\ \hline
    \textbf{RF} 
        & 582                               
        & 0.1356     
        & 582.10      
        \\ \hline
    \textbf{LR+RF}        
        & \textbf{506}                               
        & \textbf{0.1130}    
        & \textbf{506.24}   
        \\ \hline
\end{tabular}%
}
\caption{Prediction result for model 1}
\label{tab:result1}
\end{table}

\subsection{Simulation}

\begin{table*}[!t]
\centering
\resizebox{0.9\textwidth}{!}{%
\begin{tabular}{|c|c|c|c|c|c|c|c|c|c|c|c|}
\hline
\multicolumn{2}{|c|}{Game Data}        &\multicolumn{2}{c|}{Ridership}        & \multicolumn{4}{c|}{Simulated with   Actual Schedule}    & \multicolumn{4}{c|}{Simulated with   Proposed Schedule} \\ \hline
Date     & Offset & Actual & Predicted  & \# trains & Avg. WT & Std. & Avg. \% LB & \# trains & Avg. WT
& Std. & Avg. \% LB \\ \hline  \hline
6/30/18         & 0             & 7456         & 8292            & 11        & 5.7             & 3.5  & 22.3                & 12        & 3.3            & 2.8  & 1.2                 \\ \hline
7/15/18         & 25            & 8773         & 8097            & 11        & 7.2             & 4.3  & 39.2                & 12        & 4              & 3.8  & 12                  \\ \hline
9/22/18         & 0             & 7693         & 6765            & 11        & 3.8             & 2.8  & 8                   & 10        & 4              & 3.4  & 5.3                 \\ \hline
10/21/18        & 10            & 6431         & 7787            & 9         & 6.5             & 4.7  & 20.2                & 12        & 2.8            & 2    & 0                   \\ \hline
5/12/19         & 0             & 7790         & 6949            & 11        & 6.1             & 3.7  & 28.4                & 11        & 4.2            & 3.1  & 13.5                \\ \hline
7/7/19          & 25            & 7356         & 7333            & 10        & 4.7             & 3.2  & 13.6                & 11        & 3.6            & 2.9  & 3.3                 \\ \hline
8/3/19          & 0             & 6517         & 7250            & 10        & 4.4             & 3.4  & 4.8                 & 11        & 3.3            & 2.5  & 0.2                 \\ \hline
9/14/19         & 0             & 6631         & 5909            & 9         & 5.3             & 3.9  & 15.7                & 9         & 4.3            & 3.1  & 8                   \\ \hline
\end{tabular}
}
\caption{Results from two simulations for each of the eight Atlanta United Games using actual and proposed schedules assuming a max occupancy of 707 riders per vehicle.}
\label{tab:AtlUnitedSimulationTable}
\end{table*}

The forecast splits the post-game ridership in five minute bins
and can then be used to create a train schedule. The proposed schedule can be compared to the actual (recovered) schedule, giving key  insights to help dispatchers improve performance of the rail service during post-game spikes.

\paragraph{Case Study Data}
The case study is focused on the Atlanta United games with the upper deck open. Two Atlanta United games were excluded from the analysis, one because it was a playoff game and the other because there was an overlapping basketball game and no similar examples to use for the predictions. For each day, the actual schedule is recovered using the methods in Section \ref{sec:TrainLevelClustering}. Note manual adjustments are made to add in trains in the case that two trains were close together to make sure that it was a fair comparison. 
Due to the results of Section \ref{sec:train-capacity-estimation}, it is assumed that each train will fill up to a max of 707 passengers. To adjust the predictions of total ridership, from the analysis in \ref{sec:od-analysis}, it is assumed 8\% of passengers are going West, and 68\% of passengers are leaving during this peak period. Then, an extra 10\% of people are planned for the schedules to give a buffer as this is the MAPE from the previous section \ref{sec:ridership-prediction-sub}.

\paragraph{Proposed Schedule} The proposed schedule is computed by the following steps:
\begin{enumerate}
\item take the average throughput curve and divide each value by the total average throughput for this period to get the percent arrivals per bin;
\item scale down the prediction to account only for riders headed in the East direction on the rail network (92\%) and those leaving within the 120 minute time window (68\%), obtaining the predicted ridership;
\item multiply the percent arrivals per bin by the predicted total ridership;
\item design the train schedule to match this predicted demand curve by scheduling trains so that each one will depart as soon as the last passenger required to fill the train arrives.
\end{enumerate}

The proposed and actual scheduled are compared by evaluating their simulation results using Algorithm \ref{alg:left-behind-calculation} across the subset of Atlanta United game day with the upper deck seating open. The simulations focus on the peak post-game time period. It is assumed that the schedules return to normal (every 10 minutes) following this modified schedule for the peak period. 

\begin{figure}[!t]
    \centering
    \includegraphics[width=\linewidth]{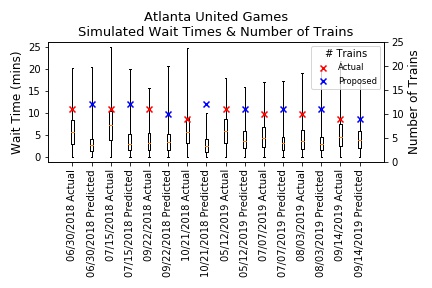}
    \caption{For each of the eight Atlanta United games with the upper deck open, the two box plots compare simulated actual train schedules vs. simulated proposed train schedules. Note the number of trains is represented with \textit{x} markers.}
    \label{fig:AtlUnitedWaitTimeBoxPolots}
\end{figure}

\begin{figure}[!t]
    \centering
    \includegraphics[width=\linewidth]{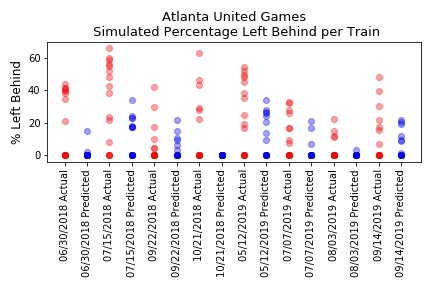}
    \caption{For each of the eight Atlanta United games with the upper deck open, the \textit{percentage of ridership left behind} for each train is plotted for both simulations. The results from the simulations of actual (recovered) schedules are in red and results from simulations with proposed schedules are in blue. Note trains near the beginning and end of the schedules tend to have 0 riders left behind.}s
    \label{fig:AtlUnitedLeftBehindPercent}
\end{figure}

\paragraph{Results}
Table \ref{tab:AtlUnitedSimulationTable} displays the results of the 16 simulations: two for each of the eight Atlanta United games of focus where the upper-decks were open. There is significant decreases to the average percent left behind (LB) on average for each train compared to the actual train, which shows the schedule is better matching the increased demand. On average about one more train is used as an extra buffer was added to protect against the times that the ridership is underpredicted. This is the number of trains that depart from DOME/GWCC within the considered time window.

Figure \ref{fig:AtlUnitedWaitTimeBoxPolots} shows boxplots of the simulated wait times using both actual and proposed schedules. The number of trains represented by the \textit{x} markers. Note that, in most cases, the maximum wait time is decreased, as well as the 75th percentile and median. 

Figure \ref{fig:AtlUnitedLeftBehindPercent} highlights the percentage of people left behind by individual trains. The proposed schedules performed significantly better in this category as the demand matches  the proposed train schedule more accurately. In the existing schedule, less than half the people at the station are able to board the train. This could lead to potential crowding and decreased customer experience as passengers have to wait multiple trains before boarding in some cases. Note trains near the beginning and end of the schedules tend to have 0 riders left behind.

\section{Discussion}

Ultimately, the planning of public transit for and around special
events is a difficult process with lots of different interacting
factors: it is almost impossible to see with perfect foresight what
will happen. However, through post-event analysis, pattern discovery,
and machine learning, public transit authorities can get closer to
that goal. By following the pipeline proposed in this paper, it
becomes possible to predict key performance metrics with reasonable
accuracy from the overall event ridership to how crowded each train
is. Much of the predictive power is derived from the behavioral
consistency of the attendees/riders: ridership is strongly correlated
to attendance and riders leave the events in highly regular patterns.

This study focused on Atlanta and sporting events and were obtained on
data for the entirety of 2018 and 2019. Over this period, there were
137 days with special events used to train and evaluate the
models. While the results achieved are of high fidelity, the accuracy can be
further improved with greater sample sizes. Mercedes-Benz Stadium and
Atlanta United did not exist before 2018, no additional historical
data is available, Similarly, Atlanta in 2020 was affected by the
COVID-19 pandemic, thus restricting the data set to only two years,
While these constraints exist in Atlanta, the methods in this paper
are valid and able to be applied elsewhere. Gathering a larger sample
size in the future, or testing these methods on other cities where the
event landscape has remained constant for longer, is another way to
improve the models. Additionally, transit agencies reaching out to event-center administers for information about ticket sales and other factors for past and future events could be very beneficial to improving attendance estimates and ridership predictions.

In the analysis, it was noted that the largest decrease in post-game ridership was to the northernmost station on the Red Line, North Springs. 
It should be noted that typically after 9:30, the Red Line service is reduced to only run part of the usual route, which typically means there is an extra transfer to make it all the way to North Springs from the nearby stadium stations. Although MARTA adds back this service during this post-game peak, it is possible that riders are not aware of this additional post-game service and the additional transfer leads to riders to look for alternative means of transportation home, such as ridesharing. It is possible that better marketing and publicity might help recapture this lost ridership, such as having an updated schedule on Google Maps for game day.

It was also noted that there is a strong correlation between station ridership and the amount of parking. It is likely that some people are choosing to park-and-ride with MARTA due to the fact that stadium parking is often expensive. Adjusting the amount of parking available at stations or the price of parking could be a interesting mechanism to help control passenger demand at various stations. Given that the majority of riders utilize the East and North directions of trains post-game. It may make sense to try to incentivize passengers to use the stations West of the stadium through increased parking or lower prices. Additionally, partnerships with ride-sharing companies in a certain areas could be another way to attract more passengers to influence the demand for ridership. 

Potential other areas for future work include building system-wide
simulations for what is likely to happen before and after games. By
using model output as simulation input, the system can be predicted
and visualized in advance of events. Other areas of future work
include finding easier ways to predict the impact of one-off or low
sample size events in the system. Consistent professional events were
easier to model than college events or national team games that
happened in Atlanta. Perhaps, analysis of examples from other cities
may be relevant in this case. It would also be interesting to expand
similar case study analysis to other types of events and other cities
to see whether the same patterns hold true. Even within Atlanta, it may be beneficial to see the impact that can be created through changes in the actual rail infrastructure such as new lines or line extensions. Yet another area to explore is the impact of post-game parties or ridesharing services on transit usage post-games.

\section*{Acknowledgments}
Many thanks to Robert Goodwin, director of Research and Analysis at
MARTA, for asking the research questions that motivated this research
and sharing MARTA's data. This research is partly supported by NSF
Leap HI award 1854684.

\bibliographystyle{IEEEtran}
\bibliography{references}

\appendices
{\onecolumn
\section{Baseline Signatures}
\begin{figure}[h!]
\begin{subfigure}[b]{0.45\textwidth}
\includegraphics[width=\linewidth]{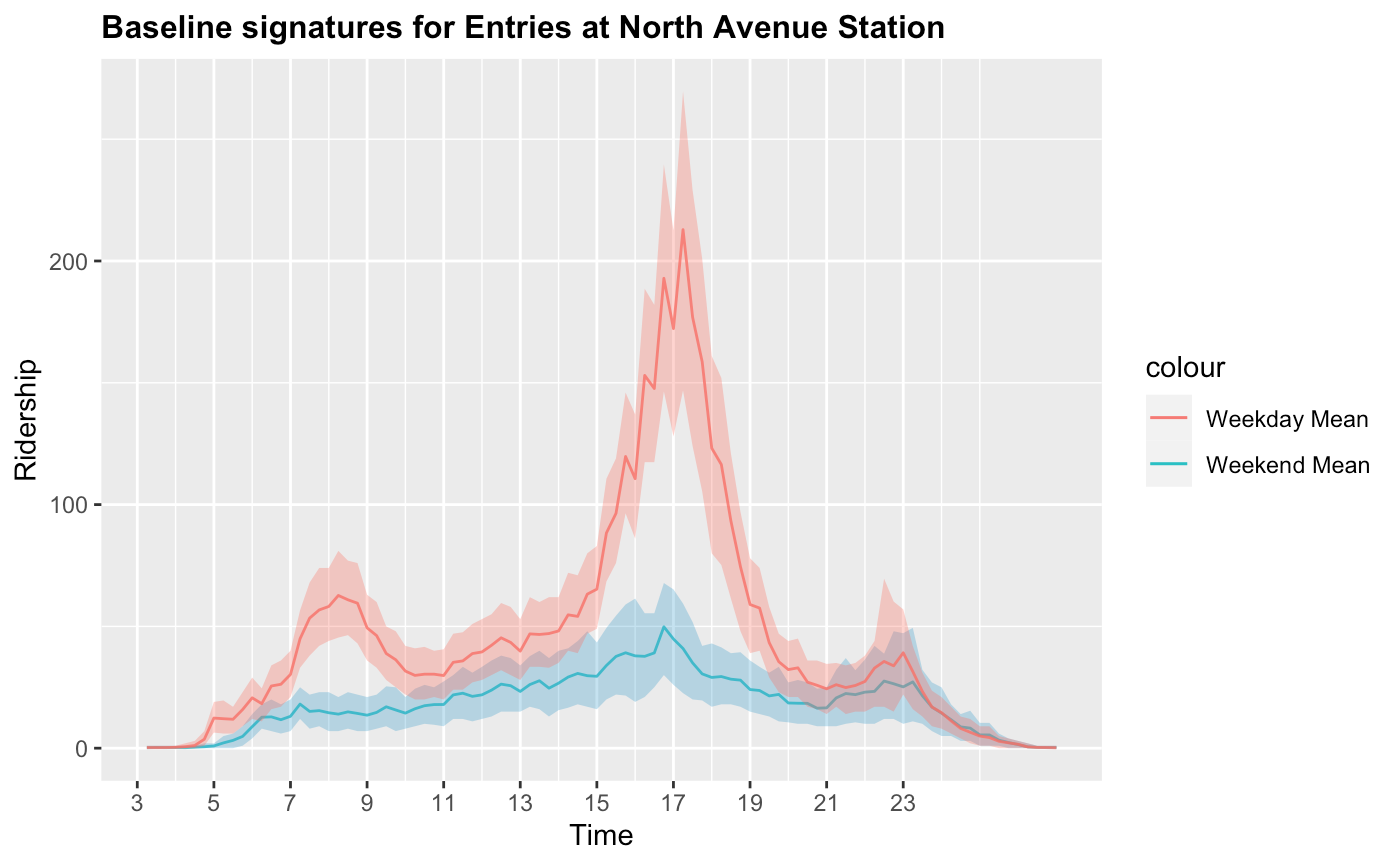}
\caption{Baseline Signature for Entries at North Ave Rail Station}
\label{fig:BaselineStart}
\end{subfigure}
\hfill
\begin{subfigure}[b]{0.45\textwidth}
\includegraphics[width=\linewidth]{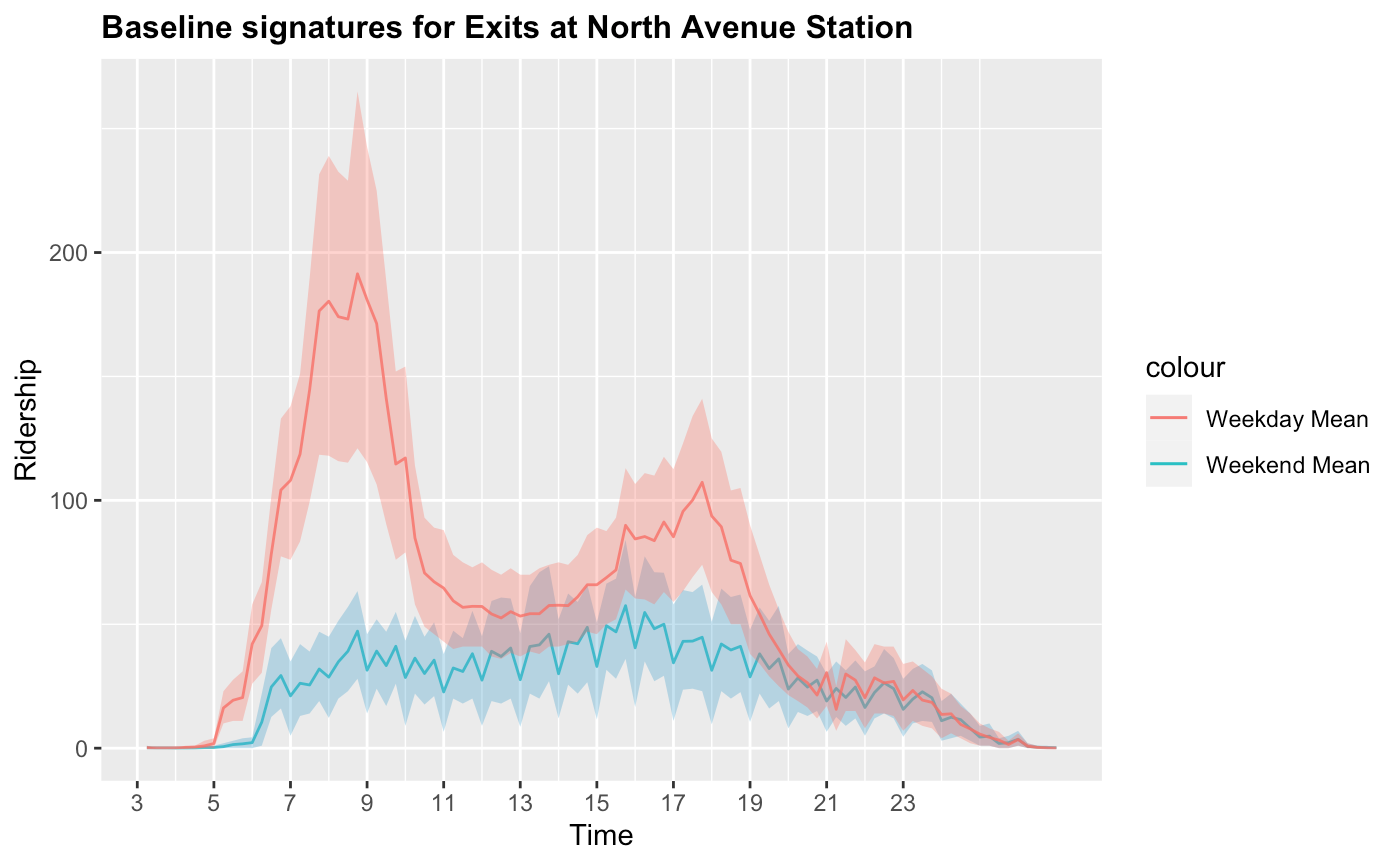}
\caption{Baseline Signature for Exits at North Ave Rail Station}
\end{subfigure}
\hfill
\begin{subfigure}[b]{0.45\textwidth}
\includegraphics[width=\linewidth]{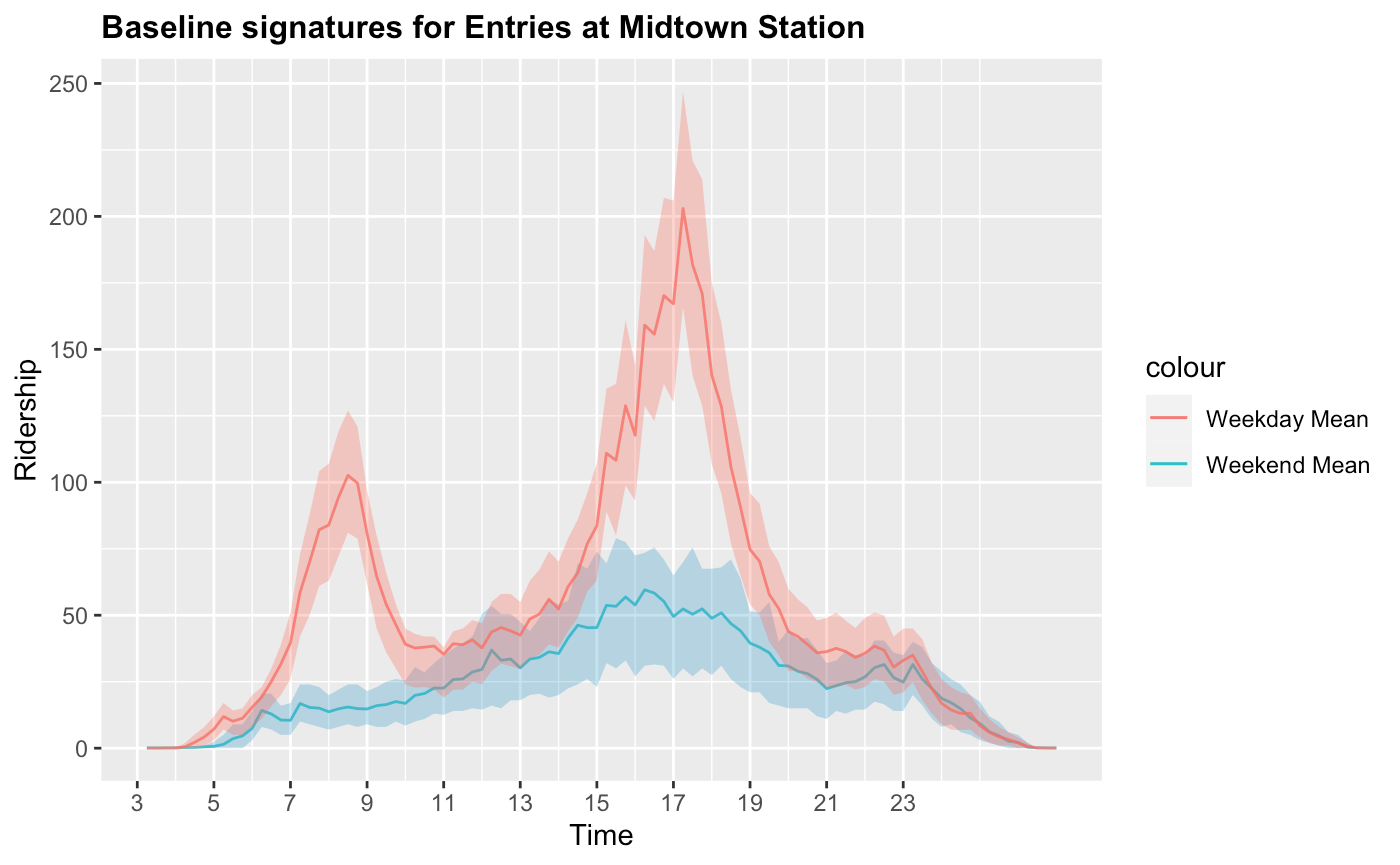}
\caption{Baseline Signature for Entries at Midtown Rail Station}
\end{subfigure}%
\hfill
\begin{subfigure}[b]{0.45\textwidth}
\includegraphics[width=\linewidth]{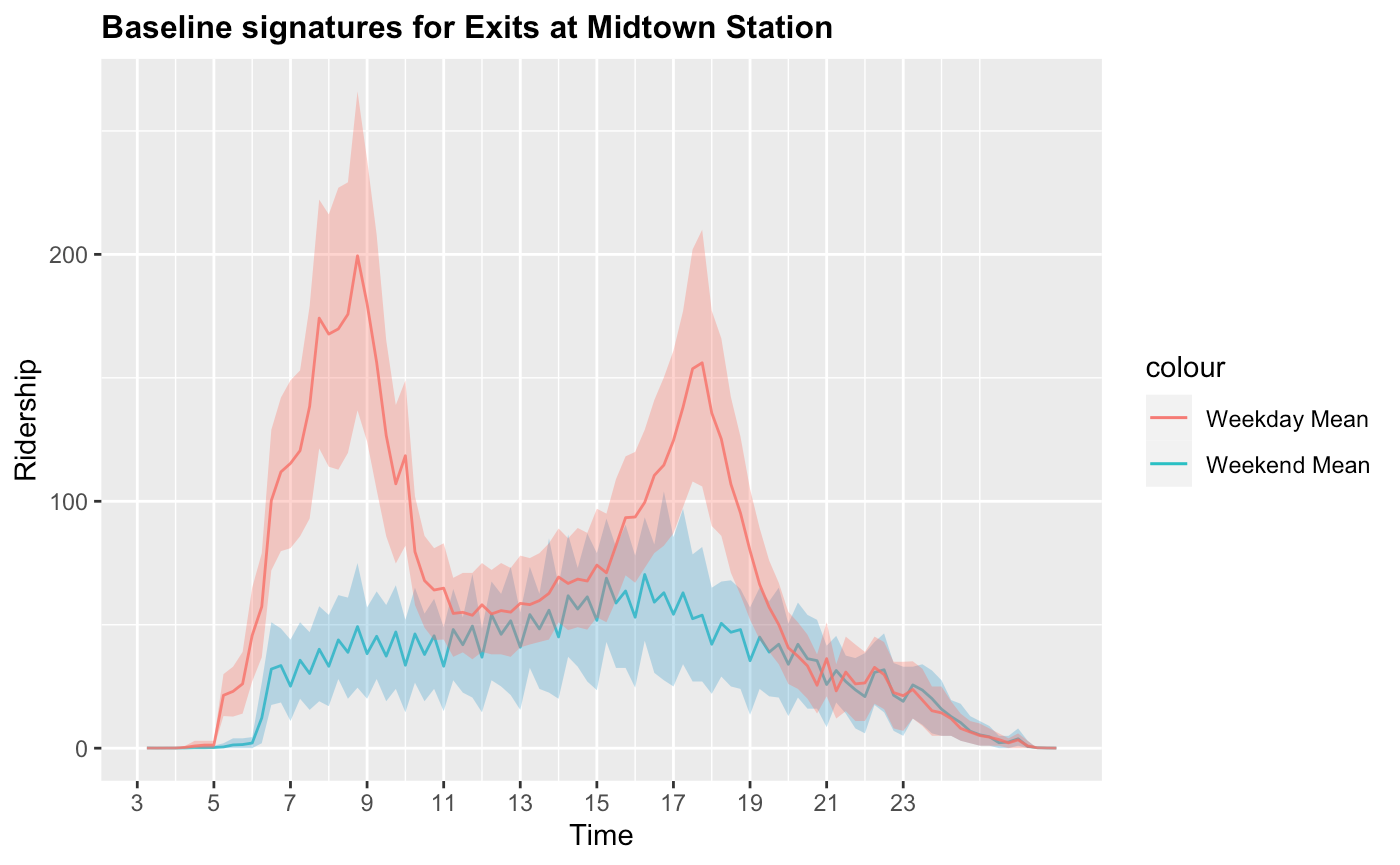}
\caption{Baseline Signature for Exits at Midtown Rail Station}
\label{fig:BaselineEnd}
\end{subfigure}
\hfill
\begin{subfigure}[b]{0.45\textwidth}
\includegraphics[width=\linewidth]{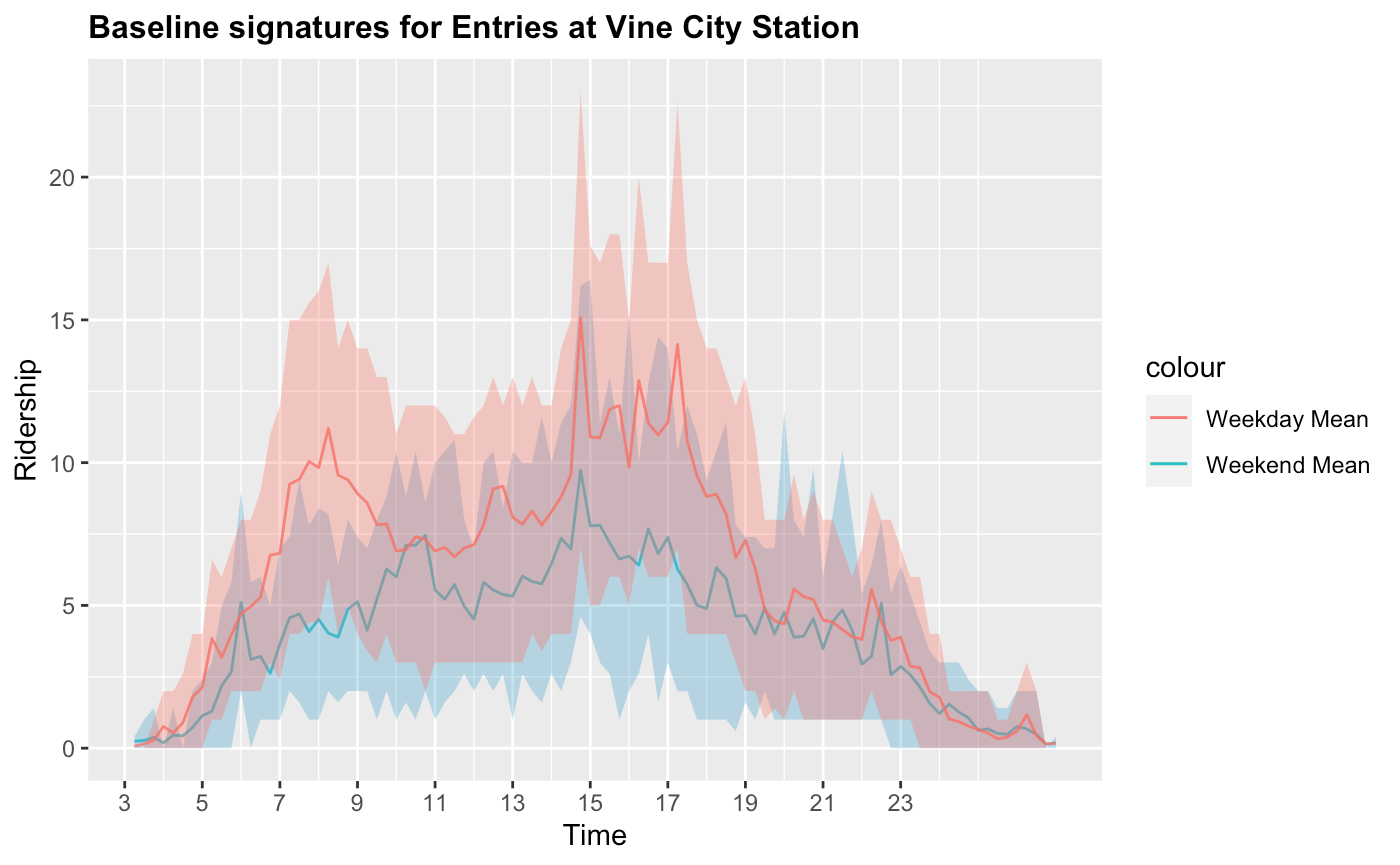}
\caption{Baseline Signature for Entries at Vine City Rail Station}
\label{fig:BaselineSignatureVineCityEntry}
\end{subfigure}%
\hfill
\begin{subfigure}[b]{0.45\textwidth}
\includegraphics[width=\linewidth]{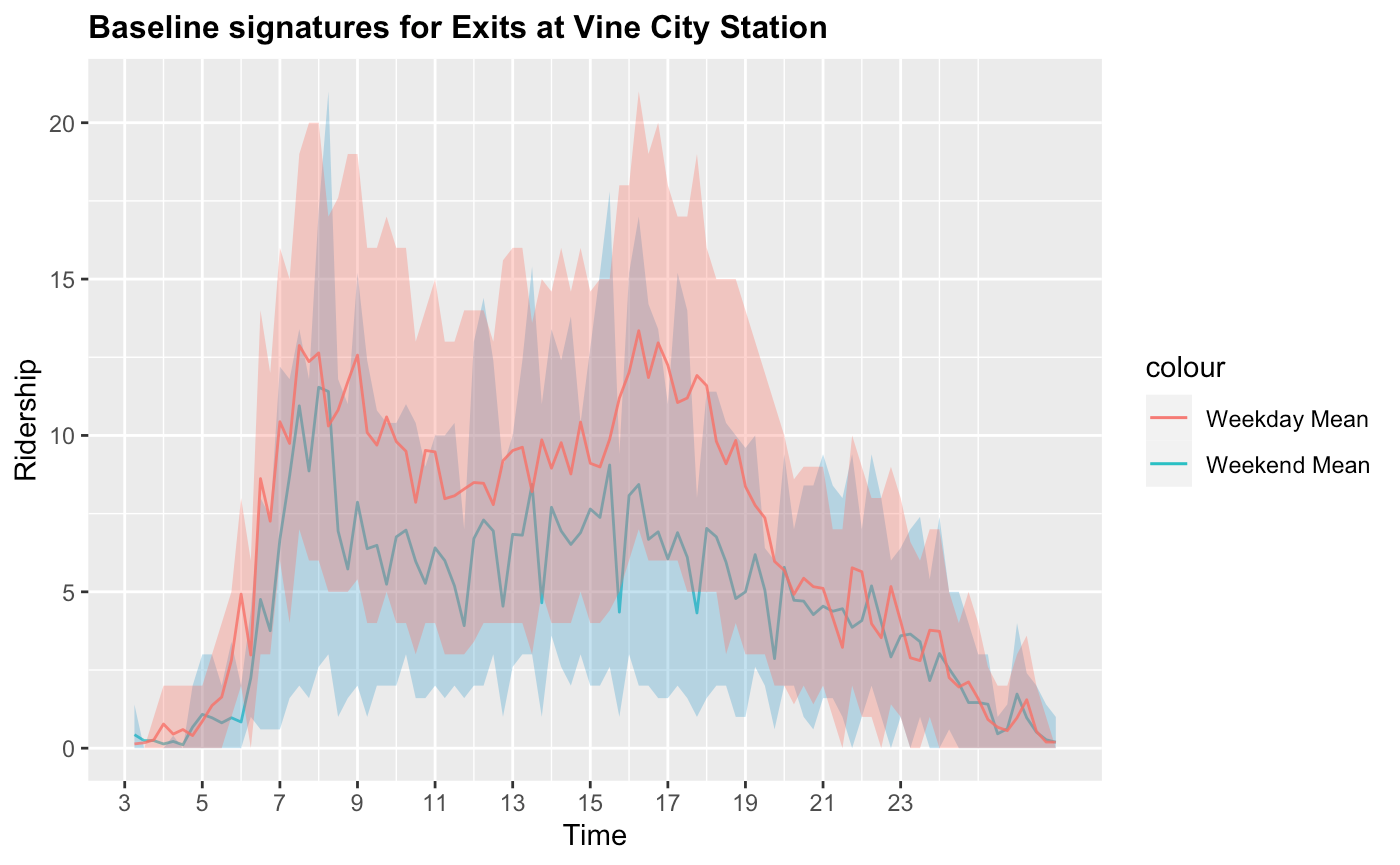}
\caption{Baseline Signature for Exits at Vine City Rail Station}
\label{fig:BaselineSignatureVineCityExit}
\end{subfigure}
\caption{Baseline (non-event days) signature graphs for three stations showing the weekday (red) and weekend (blue) signatures for entries (left) and exits (right). The shaded region represents the 10th-90th percentile range for each bin.}
\label{fig:BaselineSignatures}
\end{figure}
}

\end{document}